\documentclass[11pt]{article}
\usepackage{amsmath, upgreek}
\usepackage{amssymb}
\usepackage{amscd}
\usepackage{xspace}
\usepackage{verbatim}
	
\usepackage{color}

\newcommand{\mysection}[1]{
\section{#1}\setcounter{equation}{0}}
\title{\bf Measure data problems for a class of elliptic equations with mixed absorption-reaction}%%
\author{{\bf Marie-Fran\c{c}oise Bidaut-V\'eron\footnote{\noindent Laboratoire de Math\'{e}matiques et Physique Th\'{e}orique,
Universit\'e de Tours, 37200 Tours, France. E-mail: veronmf@univ-tours.fr},} \\{\bf Marta Garcia-Huidobro \footnote{\noindent
Departamento de Matem\'aticas, Pontifica Universidad Cat\'olica de Chile,
Casilla 307, Correo 2, Santiago de Chile. E-mail: mgarcia@mat.uc.cl}}\\
 {\bf Laurent V\'eron \footnote{\noindent
Laboratoire de Math\'{e}matiques et Physique Th\'{e}orique, Universit\'e de Tours, 37200 Tours, France. E-mail: veronl@univ-tours.fr}}\\[2mm]
}%%О©╫О©╫

\date{}
\begin{document}
 \maketitle
% \noindent{\small {\bf Abstract} We study the existence and uniqueness of  solutions of $\partial_tu-\Delta u+u^q=0$ ($q>1$) in $\Omega\times (0,\infty)$ where $\Omega\subset\mathbb R^N$ is a domain with a compact boundary, subject to the conditions $u=f\geq 0$ on $\partial\Omega\times (0,\infty)$ and the initial condition $\lim_{t\to 0}u(x,t)=\infty$. By means of Brezis' theory of maximal monotone operators in Hilbert spaces, we construct a minimal solution when $f=0$, whatever is the regularity of the boundary of the domain. When $\partial\Omega$ satisfies the parabolic Wiener criterion and $f$ is continuous, we construct a maximal solution and prove that it is the unique solution which blows-up at $t=0$.
% }

% \noindent
% {\it \footnotesize 1991 Mathematics Subject Classification}. {\scriptsize
% 35K60}.\\
% {\it \footnotesize Key words}. {\scriptsize Parabolic equations, singular solutions, semi-groups of contractions, maximal monotone operators, Wiener criterion.}
% \vspace{1mm}
% \hspace{.05in}

%% FONT commands
\newcommand{\txt}[1]{\;\text{ #1 }\;}%% Used in math only
\newcommand{\tbf}{\textbf}%% Bold face. Usage: \tbf{...}
\newcommand{\tit}{\textit}%% Italic
\newcommand{\tsc}{\textsc}%% Small caps
\newcommand{\trm}{\textrm}
\newcommand{\mbf}{\mathbf}%% Math bold
\newcommand{\mrm}{\mathrm}%% Math Roman
\newcommand{\bsym}{\boldsymbol}%% Bold math symbol
%%Macros for changing font size in math.
\newcommand{\scs}{\scriptstyle}%% as in subscript
\newcommand{\sss}{\scriptscriptstyle}%% as in sub-subscript
\newcommand{\txts}{\textstyle}
\newcommand{\dsps}{\displaystyle}
%%Macros for changing font size in text.
\newcommand{\fnz}{\footnotesize}
\newcommand{\scz}{\scriptsize}
%%\tiny<\scz<\fsz<\small<\large<\Large<\huge<\Huge
%%%%%%%%%%%%
%%%%%%%%%%%%
%% EQUATION commands
\newcommand{\be}{\begin{equation}}
\newcommand{\bel}[1]{\begin{equation}\label{#1}}
\newcommand{\ee}{\end{equation}}%% This macro does not work with amstex.
\newcommand{\eqnl}[2]{\begin{equation}\label{#1}{#2}\end{equation}}%%use not advisable; confusing
\newcommand{\barr}{\begin{eqnarray}}
\newcommand{\earr}{\end{eqnarray}}
\newcommand{\bars}{\begin{eqnarray*}}
\newcommand{\ears}{\end{eqnarray*}}
\newcommand{\nnu}{\nonumber \\}
%%%%%%%%%%%%%%%
%% Unnumbered THEOREM env.
%% New env. to be used for unnumbered theorem, lemma etc. (but with specified name)
\newtheorem{subn}{\name}
\renewcommand{\thesubn}{}
\newcommand{\bsn}[1]{\def\name{#1}\begin{subn}}
\newcommand{\esn}{\end{subn}}
%%%%%%%%%%%%%%
%% NUMBERED THEOREM env.
%% Environments: theorem, lemma, corollary defintion and related commands,
%% designed to provide consecutive numbering of these forms.
\newtheorem{sub}{\name}[section]
\newcommand{\dn}[1]{\def\name{#1}}   %used in conjuction with sub or subn.
\newcommand{\bs}{\begin{sub}}
\newcommand{\es}{\end{sub}}
\newcommand{\bsl}[1]{\begin{sub}\label{#1}}
%% the above must be preceeded by \dn (name definition),
%% however this is superceded by the list of commands bth etc.  below.
%%%%%%%%%%%%
%% NUMBERED THEOREM env. (cont.)
%% List of commands derived from 'sub' env. for theorem, lemma etc.
%% designed to provide consecutive numbering of these forms.
\newcommand{\bth}[1]{\def\name{Theorem}
\begin{sub}\label{t:#1}}
\newcommand{\blemma}[1]{\def\name{Lemma}
\begin{sub}\label{l:#1}}
\newcommand{\bcor}[1]{\def\name{Corollary}
\begin{sub}\label{c:#1}}
\newcommand{\bdef}[1]{\def\name{Definition}
\begin{sub}\label{d:#1}}
\newcommand{\bprop}[1]{\def\name{Proposition}
\begin{sub}\label{p:#1}}
%%%%%%%%%%%%%%%%%%%%%%%%%%%%%%%%%%
%% RERERENCE commands.
%% \newcommand{\R}[1]{(\ref{#1})}
\newcommand{\R}{\eqref}
\newcommand{\rth}[1]{Theorem~\ref{t:#1}}
\newcommand{\rlemma}[1]{Lemma~\ref{l:#1}}
\newcommand{\rcor}[1]{Corollary~\ref{c:#1}}
\newcommand{\rdef}[1]{Definition~\ref{d:#1}}
\newcommand{\rprop}[1]{Proposition~\ref{p:#1}}
%%%%%%%%%%%
%% ARRAY commands.
\newcommand{\BA}{\begin{array}}
\newcommand{\EA}{\end{array}}
\newcommand{\BAN}{\renewcommand{\arraystretch}{1.2}
\setlength{\arraycolsep}{2pt}\begin{array}}
\newcommand{\BAV}[2]{\renewcommand{\arraystretch}{#1}
\setlength{\arraycolsep}{#2}\begin{array}}
%Note: The first variable gives the amount of stretching: (#1) x default.
%For instance #1=1.2 means a 20% stretching. The second variable should be
%written for instance in the form  4pt ; here the default is 5pt
%\newcommand{\EAN}{\end{array}\setlength{\arraycolsep}{5pt}}
\newcommand{\BSA}{\begin{subarray}}
\newcommand{\ESA}{\end{subarray}}
%Note: These are used in subscripts as well as superscripts. They work essentially
%% like 'array'.
\newcommand{\BAL}{\begin{aligned}}
\newcommand{\EAL}{\end{aligned}}
\newcommand{\BALG}{\begin{alignat}}
\newcommand{\EALG}{\end{alignat}}%% the abbrev. does not work with latex2e
\newcommand{\BALGN}{\begin{alignat*}}
\newcommand{\EALGN}{\end{alignat*}}%% the abbrev. does not work with latex2e
%% The 'aligned' environment must be placed inside an 'equation' env.
%% in the same way as the array.
%% One could use also the 'align' env. or the 'alignat' env.
%% However in this case each line is numbered, unless '\notag' is used.
%% The 'alignat'
%% has a slightly different format (the number of columns must be specified in advance)
%% but it has the advantage that the distance between columns is at our disposition.
%% (The default would be zero distance.) Using 'alignat*' we can have the advantages
%% of alignat plus the situation where separate lines are not numbered.
%% However in this case there is no numbering at all (unless we provide a tag).
%%%%%%%%%%
%% PROOF, REMARK etc.
\newcommand{\note}[1]{\textit{#1.}\hspace{2mm}}
\newcommand{\Proof}{\note{Proof}}
\newcommand{\qeda}{\hspace{10mm}\hfill $\square$}
\newcommand{\qed}{\\
${}$ \hfill $\square$}
\newcommand{\Remark}{\note{Remark}}
%%%%%%%% Style command.
\newcommand{\modin}{$\,$\\[-4mm] \indent}
%% To be used after \mysection in order to start new line with \indent.
%%%%%%%%%%%%
%% MATHEMATICAL symbols
\newcommand{\forevery}{\quad \forall}
\newcommand{\set}[1]{\{#1\}}
\newcommand{\setdef}[2]{\{\,#1:\,#2\,\}}
\newcommand{\setm}[2]{\{\,#1\mid #2\,\}}
%% Arrows
\newcommand{\mt}{\mapsto}
\newcommand{\lra}{\longrightarrow}
\newcommand{\lla}{\longleftarrow}
\newcommand{\llra}{\longleftrightarrow}
\newcommand{\Lra}{\Longrightarrow}
\newcommand{\Lla}{\Longleftarrow}
\newcommand{\Llra}{\Longleftrightarrow}
\newcommand{\warrow}{\rightharpoonup}
%% Brackets, delimiters
\newcommand{
\paran}[1]{\left (#1 \right )}%% adjustable parantheses
\newcommand{\sqbr}[1]{\left [#1 \right ]}%% adjustable square brackets
\newcommand{\curlybr}[1]{\left \{#1 \right \}}%% adjustable curly brackets
\newcommand{\abs}[1]{\left |#1\right |}%% adjustable vertical delimiters
\newcommand{\norm}[1]{\left \|#1\right \|}%% adjustable norm
\newcommand{
\paranb}[1]{\big (#1 \big )}%% non-adjustable parantheses (big)
\newcommand{\lsqbrb}[1]{\big [#1 \big ]}%% non-adjustable square brackets (big)
\newcommand{\lcurlybrb}[1]{\big \{#1 \big \}}%% non-adjustable curly brackets (big)
\newcommand{\absb}[1]{\big |#1\big |}%% non-adjustable vertical delimiters (big)
\newcommand{\normb}[1]{\big \|#1\big \|}%% non-adjustable norm (big)
\newcommand{
\paranB}[1]{\Big (#1 \Big )}%% non-adjustable parantheses (Big)
\newcommand{\absB}[1]{\Big |#1\Big |}%% non-adjustable vertical delimiters (Big)
\newcommand{\normB}[1]{\Big \|#1\Big \|}%% non-adjustable norm (Big)
\newcommand{\produal}[1]{\langle #1 \rangle}%% the pairing of X' and X 

%%%%%%%%%%%%%%%%%
%% Adjustable parantheses etc. in a different DEFINITION format.
%\def\adp(#1){\left (#1 \right )}%% adjustable parantheses
%\def\adsb(#1){\left [#1\right ]}%% adjustable square brackets
%\def\adcb(#1){\left \{#1\right \}}%% adjustable curly brackets
%\def\abs|#1|{\left |#1\right |}%% adjustable vertical delimiters
%%%%%%%%%%%%%%%%
%% More mathematical symbols
\newcommand{\thkl}{\rule[-.5mm]{.3mm}{3mm}}
\newcommand{\thknorm}[1]{\thkl #1 \thkl\,}
\newcommand{\trinorm}[1]{|\!|\!| #1 |\!|\!|\,}
\newcommand{\bang}[1]{\langle #1 \rangle}%% angle bracket
\def\angb<#1>{\langle #1 \rangle}%% angle bracket
%% The two last lines yield the same result.
%% The second is used as follows: \angb<a,b>
\newcommand{\vstrut}[1]{\rule{0mm}{#1}}
\newcommand{\rec}[1]{\frac{1}{#1}}
%% OPERATOR names.
%% OPERATOR names.
\newcommand{\opname}[1]{\mbox{\rm #1}\,}
\newcommand{\supp}{\opname{supp}}
\newcommand{\dist}{\opname{dist}}
\newcommand{\myfrac}[2]{{\displaystyle \frac{#1}{#2} }}
\newcommand{\avint}[2]{{\displaystyle-\!\!\!\!\myint{#1}{#2}\!\!\!\!\!\!\!\!\!\!\!\!\!\!\!-\,\,\;\;\;\;}}
\newcommand{\myint}[2]{{\displaystyle \int_{#1}^{#2}}}
\newcommand{\mysum}[2]{{\displaystyle \sum_{#1}^{#2}}}
\newcommand {\dint}{{\displaystyle \myint\!\!\myint}}%%%%%%%%%%
%%%%%%% SPACE commands
\newcommand{\q}{\quad}
\newcommand{\qq}{\qquad}
\newcommand{\hsp}[1]{\hspace{#1mm}}
\newcommand{\vsp}[1]{\vspace{#1mm}}
%%%%%%%%%%%
%% ABREVIATIONS
\newcommand{\ity}{\infty}
\newcommand{\prt}{\partial}
\newcommand{\sms}{\setminus}
\newcommand{\ems}{\emptyset}
\newcommand{\ti}{\times}
\newcommand{\pr}{^\prime}
\newcommand{\ppr}{^{\prime\prime}}
\newcommand{\tl}{\tilde}
\newcommand{\sbs}{\subset}
\newcommand{\sbeq}{\subseteq}
\newcommand{\nind}{\noindent}
\newcommand{\ind}{\indent}
\newcommand{\ovl}{\overline}
\newcommand{\unl}{\underline}
\newcommand{\nin}{\not\in}
\newcommand{\pfrac}[2]{\genfrac{(}{)}{}{}{#1}{#2}}% frac with parantheses.
%%%%%%%%%%%
%%%%%%%%%%%%%

%%Macros for Greek letters.
\def\ga{\alpha}     \def\gb{\beta}       \def\gg{\gamma}
\def\gc{\chi}       \def\gd{\delta}      \def\ge{\epsilon}
\def\gth{\theta}                         \def\vge{\varepsilon}
\def\gf{\phi}       \def\vgf{\varphi}    \def\gh{\eta}
\def\gi{\iota}      \def\gk{\kappa}      \def\gl{\lambda}
\def\gm{\mu}        \def\gn{\nu}         \def\gp{\pi}
\def\vgp{\varpi}    \def\gr{\rho}        \def\vgr{\varrho}
\def\gs{\sigma}     \def\vgs{\varsigma}  \def\gt{\tau}
\def\gu{\upsilon}   \def\gv{\vartheta}   \def\gw{\omega}
\def\gx{\xi}        \def\gy{\psi}        \def\gz{\zeta}
\def\Gg{\Gamma}     \def\Gd{\Delta}      \def\Gf{\Phi}
\def\Gth{\Theta}
\def\Gl{\Lambda}    \def\Gs{\Sigma}      \def\Gp{\Pi}
\def\Gw{\Omega}     \def\Gx{\Xi}         \def\Gy{\Psi}

%%Macros for calligraphic letters.
\def\CS{{\mathcal S}}   \def\CM{{\mathcal M}}   \def\CN{{\mathcal N}}
\def\CR{{\mathcal R}}   \def\CO{{\mathcal O}}   \def\CP{{\mathcal P}}
\def\CA{{\mathcal A}}   \def\CB{{\mathcal B}}   \def\CC{{\mathcal C}}
\def\CD{{\mathcal D}}   \def\CE{{\mathcal E}}   \def\CF{{\mathcal F}}
\def\CG{{\mathcal G}}   \def\CH{{\mathcal H}}   \def\CI{{\mathcal I}}
\def\CJ{{\mathcal J}}   \def\CK{{\mathcal K}}   \def\CL{{\mathcal L}}
\def\CT{{\mathcal T}}   \def\CU{{\mathcal U}}   \def\CV{{\mathcal V}}
\def\CZ{{\mathcal Z}}   \def\CX{{\mathcal X}}   \def\CY{{\mathcal Y}}
\def\CW{{\mathcal W}} \def\CQ{{\mathcal Q}}
%%%%%
%%Macros for 'blackboard' letters (See (27) for display.)
\def\BBA {\mathbb A}   \def\BBb {\mathbb B}    \def\BBC {\mathbb C}
\def\BBD {\mathbb D}   \def\BBE {\mathbb E}    \def\BBF {\mathbb F}
\def\BBG {\mathbb G}   \def\BBH {\mathbb H}    \def\BBI {\mathbb I}
\def\BBJ {\mathbb J}   \def\BBK {\mathbb K}    \def\BBL {\mathbb L}
\def\BBM {\mathbb M}   \def\BBN {\mathbb N}    \def\BBO {\mathbb O}
\def\BBP {\mathbb P}   \def\BBR {\mathbb R}    \def\BBS {\mathbb S}
\def\BBT {\mathbb T}   \def\BBU {\mathbb U}    \def\BBV {\mathbb V}
\def\BBW {\mathbb W}   \def\BBX {\mathbb X}    \def\BBY {\mathbb Y}
\def\BBZ {\mathbb Z}   \def\BBQ {\mathbb Q}

%%Macros for Ghotic (Fraktur) letters.
\def\GTA {\mathfrak A}   \def\GTB {\mathfrak B}    \def\GTC {\mathfrak C}
\def\GTD {\mathfrak D}   \def\GTE {\mathfrak E}    \def\GTF {\mathfrak F}
\def\GTG {\mathfrak G}   \def\GTH {\mathfrak H}    \def\GTI {\mathfrak I}
\def\GTJ {\mathfrak J}   \def\GTK {\mathfrak K}    \def\GTL {\mathfrak L}
\def\GTM {\mathfrak M}   \def\GTN {\mathfrak N}    \def\GTO {\mathfrak O}
\def\GTP {\mathfrak P}   \def\GTR {\mathfrak R}    \def\GTS {\mathfrak S}
\def\GTT {\mathfrak T}   \def\GTU {\mathfrak U}    \def\GTV {\mathfrak V}
\def\GTW {\mathfrak W}   \def\GTX {\mathfrak X}    \def\GTY {\mathfrak Y}
\def\GTZ {\mathfrak Z}   \def\GTQ {\mathfrak Q}

\font\Sym= msam10 % special symbols
\def\SYM#1{\hbox{\Sym #1}}
\newcommand{\bdw}{\prt\Gw\xspace}
\maketitle\medskip

{\abstract We study the existence of nonnegative solutions to the Dirichlet problem $\CL^{_{^M}}_{p,q}u:=-\Gd u+u^p-M|\nabla u|^q=\gm$ in a domain $\Gw\subset\BBR^N$ 
where $\gm$ is a nonnegative Radon measure, when $p>1$, $q>1$ and $M\geq 0$. We also give conditions under which nonnegative solutions
of $\CL^{_{^M}}_{p,q}u=0$ in $\Gw\setminus K$ where $K$ is a compact subset of $\Gw$ can be extended as a solution of the same equation in 
$\Gw$. 

\nind {\it 2010 Mathematics Subject Classification:} 35J62-35J66-31C15-28A12\smallskip

\nind{\it  Keywords:} Elliptic equations, singularities, Bessel capacities, Riesz potential, maximal functions.

\tableofcontents
\date{}

%%%%%%%%%%%%%%%%%%%%%%%%%%%%%%%%%%%%%%%%%%%%%%%%%%%%%%%%%%%%%%%%%%%%%%%%%%%%%%%%%%%%%%%%%%%%%%%%%%%%%%%%%%%%%%%%%%%%%%%%%%%%%%%%%%%%%%%%%%%%%%%%%%%%%%%%%%%%%%%%%%%%%%%%%%%%%%%%%%%%
%%%%%%%%%%%%%%%%%%%%%%%%%%%%%%%%%%%%%%%%%%%%%%%%%%%%%%%%%%%%%SECTION--INTRODUCTION%%%%%%%%%%%%%%%%%%%%%%%%%%%%%%%%%%%%%%%%%%%%%%%%%%%%%%%%%%%%%%%%%%%%%%%%%%%%%%%%%%%%%%%%%%%%%%%%%%%%%%%%%%%%%%%%%%%%%%%%%%%%%%%%%%%%%%%%%%%%%%%%%%%%%%%%%%%%%%%%%%%%

%%%%%%%%%%%%%%%%%%%%%%%%%%%%%%%%%%%%%%%%%%%%%%%%%%%%%%%%%%
%%%%%%%%%%%%%%%%%%%%%%%%%%%%%%%%%%%%%%%%%%%%%%%%%%%%%%%%%%%%%%%%%%%%%%%%%%%%%%%%%%%%%%%%%%%%%%%%%%%%%%%%%%%%%%%%%%%%
%%%%%%%%%%%%%%%%%%%%%%%%%%%%%%%%%%%%%%%%%%%%%%%%%%%%%%%%%%%%%%%%%%%%%%%%%%%%%%%%%%%%%%%%%%%%%%%%%%%%%%%%%%%%%%%%%%%%%%%%%%%%%%%%%%%%%%%%%%%%%%%%%%%%%%%%%%%%%%%%%%%%%%%%%%%%%%%%%%%%%%%%%%%%%%%%%%%%%%%%%%%%%%%%%%%%%%%%%%%%%%%%%%%%%%%%%%%%%%%%%%%%%%%%%%%%%%%%%%%%%%%%%%%%%%%%%%%%%%%%%%%%%%%%%%%%%%%%%%%%%%%%%%%%%%%%%%%%%%%%

%%%%%%%%%%%%%%%%%%%%%%%%%%%%%%%%%%%%%%%%%%%%%%%%%%%%%%%%%%%%%%%%%%%%%%%%%%%%%%%%%%%%%%%%%%%%%%%%%%%%%%%%%%%%%%%%%%%%%%%%%%%%%%%%%%%%%%%%%%%%%%%%%%%%%%%%%%%%%%%%%%%%%%%%%%%%%%%%%%%%%%%%%%%%%%%%%%%%%%%%%%%%%%%%%%%%%%%%%%%%%%%%%%%%%%%%%%%%%%%%%%%%
%%%%%%%%%%%%%%%%%%%%%%%%%%%%%%%%%%%%%%%%%%%%%%%%%%%%%%%%%%%%%%%%%%%%%%%%%%%%%%%%%%%%%%%%%%%%%%%%%%%%%%%%%%%%%%%%%%%%
\mysection{Introduction}
Let $\Gw$ be a bounded domain  of $\BBR^N$, $N\geq 2$, and $\CL^{_{^M}}_{p,q}$ be the operator 
\begin{equation}\label{Z1}
\BA{lll}
u\mapsto \CL^{_{^M}}_{p,q} u:=-\Gd u+|u|^{p-1}u-M|\nabla u|^{q}\quad\text{for all }u\in C^2(\Gw)
\EA
\end{equation}
where $M\geq 0$ and $p,q>1$. We first provide an {\it a priori} estimate for a positive solution of $(\ref{Z1})$ and its gradient in the range $1<q<p$. Then we
study under what conditions on the parameters  any solution of 
\begin{equation}\label{Z2-0}
\BA{lll}
\CL^{_{^M}}_{p,q} u=0\qquad\text{in }\Gw\setminus K,
\EA
\ee
where $K$ is a compact subset of $\Gw$, can be extended as a solution of the same equation in whole $\Gw$, and if it is the case, whether 
the solution is bounded or not in $\Gw$. We also consider the Dirichlet problem with measure data
\begin{equation}\label{Z2}
\BA{lll}
\CL^{_{^M}}_{p,q} u=\gm\qquad &\text{in }\Gw\\
\phantom{\CL^{_{^M}}_{p,q}}u=0\qquad &\text{in }\prt\Gw,
\EA
\end{equation}
where $\gm$ is a nonnegative bounded Radon measure in $\Gw$ and exhibit conditions which guarantee the existence of nonnegative solutions to this problem.\smallskip

If $M=0$, $\CL^{_{^M}}_{p,q}$ reduces to the Emden-Fowler operator 
\begin{equation}\label{Z3}
\BA{lll}
u\mapsto \CL_{p}u:=-\Gd u+|u|^{p-1}u.
\EA
\end{equation}
Singularity problems for solutions of $\CL_{p}u=0$ have been investigated since fourty years, starting with the work of Brezis and V\'eron \cite{BrVe} who gave conditions for the removability of an isolated singularity. Later on Baras and Pierre \cite{BaPi} extended  the result in \cite{BrVe} to more general removable sets, introducing the good framework. They obtained a necessary and sufficient condition expressed in terms of the Bessel capacities
$cap^{\BBR^N}_{2,p'}$ ($p'=\frac p{p-1}$) both for the removability of compact subsets of $\Gw$ and the solvability of the associated Dirichlet problem with measure data. \
Another class of operator strongly related to $\CL^{_{^M}}_{p,q} $ is the Riccati operator 
\begin{equation}\label{Z4}
\BA{lll}
u\mapsto \CR^{_{^M}}_{q} u:=-\Gd u-M|\nabla u|^{q}.
\EA
\end{equation}
The Dirichlet problem with measure data 
\begin{equation}\label{Z5}
\BA{lll}
-\Gd u-M|\nabla u|^{q}=\gm&\qquad\text{in }\Gw\\
\phantom{-\Gd-M|\nabla u|^{q}}
u=0&\qquad\text{on }\prt\Gw
\EA
\end{equation}
has been studied by Maz'ya and Verbitsky \cite{MaVer} and Hansson, Maz'ya and Verbitsky \cite{HMV} when  $q>2$ (and also in $\BBR^N$ when $q>1$) and Phuc \cite{Ph}. Their results necessitate an extensive use of Riesz potentials. 

When $M>0$  there is a balance between the absorption term $|u|^{p-1}u$ and the source term $M|\nabla u|^{q}$, and this interaction is the origin of many unexpected new effects. In the study of singularity problems the effect of the diffusion can be neglectable compared to the 
 two nonlinear terms. The scale of the two opposed reaction terms depends upon the position of $q$ with respect to $\frac{2p}{p+1}$. This is due to the fact that if $q=\frac{2p}{p+1}$,  $(\ref{Z1})$ is equivariant with respect to the scaling transformation $T_\ell$ defined for $\ell>0$ by 
 \begin{equation}\label{Z6}
\BA{lll}
T_\ell[u](x)=\ell^{\frac{2}{p-1}}u(\ell x)=\ell^{\frac{2-q}{q-1}}u(\ell x).
\EA
\end{equation}
\nind If $q<\frac{2p}{p+1}$, the absorption term is dominant and the behaviour of the singular solutions is modelled by the equation 
$\CL_{p}u=0$ studied in \cite{Vesingsol}. 
If $q>\frac{2p}{p+1}$ , the diffusion is negligible and the behaviour of the singular solutions is modelled by positive separable solutions of
$\CE^{_{^M}}_{p,q}u=0$ where $\CE^{_{^M}}_{p,q}$ is an eikonal type operator defined by
 \begin{equation}\label{Z6-1}
\BA{lll}
\CE^{_{^M}}_{p,q}u=u^p-M|\nabla u|^q.
\EA
\end{equation}
This problem is studied in the forthcoming article \cite{BVGHV3}. If $q=\frac{2p}{p+1}$, the coefficient $M>0$ plays a fundamental role in the properties of the set of solutions, in particular for the existence of singular solutions and removable singularities; this is not the case when $q\neq\frac{2p}{p+1}$ since by an homothety $M$ can be assumed to be equal to $1$. \smallskip

Brezis and V\'eron proved in \cite{BrVe} that isolated singularities of solutions of $\CL_{p}u=0$ are removable when $p\geq \frac N{N-2}$. The removability property has been extended to more general sets using a capacity framework in \cite{BaPi}. Using a change of variable inspired by \cite{BVGHV1} where boundary singularities of solutions of $(\ref{Z6-1})$ are considered we prove a series of removability results for solutions of 
 \begin{equation}\label{Z11}\CL_{p,q}^{_{^M}}u=0.
 \end{equation}

\bth{Remov-1} Assume $0\in\Gw\subset\BBR^N$, $N\geq 3$, $M>0$, 
 $p\geq \frac N{N-2}$, $1<q\leq \frac {2p}{p+1}$ and $(p,q)\neq(\frac N{N-2},\frac {N}{N-1} )$.
Then any nonnegative solution $u\in C^2(\Gw\setminus\{0\})$ of $(\ref{Z11})$ in $\Gw\setminus\{0\}$ belongs to 
$W^{1,q}_{loc}(\Gw)\cap L^p_{loc}(\Gw)$, and it can be extended as a weak solution of $(\ref{Z11})$ in $\Gw$.\smallskip

\nind Furthermore, if we assume either\smallskip

\nind (i) $p\geq\frac N{N-2}$ and $1<q< \frac {2p}{p+1}$, or \smallskip

\nind (ii) $p>\frac N{N-2}$, $q= \frac {2p}{p+1}$ and
 \begin{equation}\label{Z12}\BA{lll}
M<m^*:=(p+1)\left(\myfrac{(N-2)p-N}{2p}\right)^{\frac{p}{p+1}},
\EA
\end{equation}
then $u\in C^2(\Gw)$.
\es

The existence of radial singular solutions when $(p,q)=(\frac N{N-2},\frac {N}{N-1} )$ and $M>0$, or when $p>\frac N{N-2}$, $q= \frac {2p}{p+1}$ and $M\geq m^*$ shows the optimality of the statements (see \cite{BVGHV3}). A series of pointwise {\it a priori} estimates concerning $u$ and $\nabla u$ are presented in the first section. They are obtained by a combination of Keller-Osserman type estimates, rescaling techniques and Bernstein method. They play a key role  for analyzing the case $p=\frac N{N-2}$ in the previous theorem, and will be of fundamental importance in the forthcoming paper \cite{BVGHV3}.
\smallskip

The method introduced in the proof of \rth{Remov-1} combined with the result of \cite{BaPi} yields a more general removability result. For such a task we denote by $cap^{\BBR^N}_{k,b}$ the Bessel capacity relative to $\BBR^N$ with order $k>0$ and power $b\in (1,\infty)$. If $k\in\BBN^*$  it coincides with the Sobolev capacity associated to the space $W^{k,b}(\BBR^N)$ by Calderon's theorem (see e.g. \cite{AdHe} for a detailed presentation).

\bth{Remov-2} Let $p>\frac{N}{N-2}$ and $\frac{N}{N-2}<r<p$. Suppose that one of the following conditions is verified:\smallskip

\nind (i) either $q=\frac{2p}{p+1}$ and 
 \begin{equation}\label{Z13}\BA{lll}
0<M<m^*(r):=(p+1)\left(\myfrac{p-r}{p(r-1)}\right)^{\frac{p}{p+1}},
\EA
\end{equation}
\smallskip

\nind (ii) or $1<q<\frac{2p}{p+1}$ and $M>0$.\smallskip

\nind Then, if $K$ is a compact subset of $\Gw$ such that $cap^{\BBR^N}_{2,r'}(K)=0$, any nonnegative solution $u\in C^2(\Gw\setminus K)$ of 
$(\ref{Z11})$ in $\Gw\setminus K$ can be extended to $\Gw$ as a solution still denoted by $u$ in the sense of distributions in $\Gw$. Furthermore, if $r\leq \frac{2N}{N-2}$, then $u\in C^2(\Gw)$.
\es

Next we obtain sufficient conditions on a positive measure in $\Gw$  in order $(\ref{Z2})$ be solvable. In the sequel we assume that  $\Gw\subset\BBR^N$, $N\geq2$, is a bounded smooth domain. We denote by $\mathfrak M(\Gw)$ (resp. $\mathfrak M^b(\Gw)$) the set of Radon measures  (resp. bounded Radon measures) in $\Gw$ and by $\mathfrak M_+(\Gw)$ (resp. $\mathfrak M_+^b(\Gw)$) its positive cone. The total variation norm of a bounded measure $\gm$ is $\norm\gm_{\mathfrak M}$.

Since for any $\gm\in \mathfrak M^b_+(\Gw)$ the nonnegative solutions of $\CL_pv=\gm$ and  $\CR^M_qw=\gm$ are respectively a subsolution and a supersolution of equation $(\ref{Z2})$ and they satisfy $0\leq v\leq w$, the construction of $v$ and $w$ is the key-stone for solving  $(\ref{Z2})$. It is known that these two problems
can be solved when the measure $\gm$ satisfies some continuity properties with respect to some specific capacities.\smallskip

\bth{meas-1} Assume $p>1$, $1<q<2$. Let $\gm\in\mathfrak M_+^b(\Gw) $. If $\gm$ satisfies
 \begin{equation}\label{Z7}
\BA{lll}
\gm(E)\leq C\min\left\{cap^{\BBR^N}_{2,p'}(E),cap^{\BBR^N}_{1,q'}(E)\right\} \quad\text{for all Borel sets }E\subset \Gw,
\EA
\end{equation}
 there is a constant $c_0>0$ such that for any $0\leq c\leq c_0$ there exists a function $u \in W^{1,q}_0(\Gw)\cap  L^p(\Gw)$, $u\geq 0$,  satisfying
 \begin{equation}\label{Z8}
\BA{lll}
-\myint{\Gw}{}u\Gd\gz dx+\myint{\Gw}{}\left(u^p-M|\nabla u|^q \right)\gz dx=c\myint{\Gw}{}\gz d\gm\quad\text{for all }\gz\in C_c^2(\overline\Gw).
\EA
\end{equation}
\es

The condition on the measure is satisfied if $W^{-1,q}(\Gw)\hookrightarrow W^{-2,p}(\Gw)$, and we prove the following:

\bcor{meas-cor1} Let $\frac{Np}{N+p}\leq q<2$ and $\gm\in\mathfrak M_+^b(\Gw) $ be such that 
 \begin{equation}\label{Z9-1}\BA{lll}
 \gm(E)\leq Ccap^{\BBR^N}_{1,q'}(E) \quad\text{for all Borel set }E\subset \Gw,
\EA
\end{equation}
for some $C>0$, then there exists a constant $c_1>0$ such that for any $0\leq c\leq c_1 $ there exists a nonnegative function $u \in W^{1,q}_0(\Gw)\cap  L^p(\Gw)$  satisfying $(\ref{Z8})$. 
\es

By comparison results between capacities we have another type of result:

\bcor{meas-2} Let $\frac{N}{N-1}\leq q\leq \frac{2p}{p+1}$. If $\gm\in\mathfrak M_+^b(\Gw) $ satisfies,
 \begin{equation}\label{Z9-2}\BA{lll}
 \gm(E)\leq Ccap^{\BBR^N}_{2,p'}(E)\quad\text{for all Borel set }E\subset \Gw,
\EA
\end{equation}
for some $C>0$, then there exists $c_2>0$ such that for any $0\leq c\leq c_2 $ there exists a nonnegative function $u \in W^{1,q}_0(\Gw)\cap  L^p(\Gw)$ satisfying $(\ref{Z8})$.
\es

As an application of the previous results, we prove the following

\bcor{meas-3} Let $p>1$, $1< q<2$ and $\gm\in\mathfrak M_+^b(\Gw)$. There exists a function $u \in W^{1,q}_0(\Gw)\cap  L^p(\Gw)$ solution of $(\ref{Z8})$ if one of the following conditions is satisfied:\smallskip

\nind (i) When $p<\frac{N}{N-2}$ and $q<\frac{N}{N-1}$, if $\norm\gm_{\mathfrak M}\leq c_3$ for some $c_3>0$.
\smallskip

\nind (ii) When $p<\frac{N}{N-2}$ and $\frac{N}{N-1}\leq q<2$, if $\gm$ satisfies $(\ref{Z9-1})$ for some $C>0$. In that case  there exists $c_4>0$ such that there must hold $0<c<c_4$
in problem $(\ref{Z7})$.
\smallskip

\nind (iii) When $p\geq \frac{N}{N-2}$ and $q<\frac{N}{N-1}$, if $\norm\gm_{\mathfrak M}\leq c^*_4M^{-\frac{1}{q-1}}$ for some $c^*_4=c^*_4(N,q,\Gw)>0$ which can be estimated, and if 
 \begin{equation}\label{Z10}\BA{lll}
 \gm(E)=0\quad\text{for all Borel set $\subset\Gw$ such that }\, cap^{\BBR^N}_{2,p'}(E)=0.
\EA
\end{equation}
\es

In the case (i) we show in a forthcoming article \cite{BVGHV3} and by a completely different method that there is no restriction on $c$ if $\gm=c\gd_a$ for some $a\in\Gw$. In the above mentioned article we construct many types of local or global singular solutions using methods inherited from dynamical systems.\medskip

\noindent{\bf Acknowledgements} This article has been prepared with the support of the  FONDECYT grants 1210241 and 1190102 for the three authors.
%\noindent{\bf Acknowledgements} This article has been prepared with the support of the  FONDECYT grants 1160540 and 1190102 for the three authors.
%%%%%%%%%%%%%%%%%%%%%%%%%%%%%%%%%%%%%%%%%%%%%%%%%%%%%%%%%%%%%%%%%%%%%%%%%%%%%%%%%%%%%%%%%%%%%%%%%%%%%%%%%%%%%%%%%%%%%%%%%%%%%%%%%%%%%%%%%%%%%%%%%%%%%%%%%%%%%%%%%%%%%%%%%%%%%%%%%%%%%%%%%%%%%%%%%%%%%%%%%%%%%%%%%%%%%%%%
%%%%%%%%%%%%%%%%%%%%%%%%%%%%%%%%%%%%%%%%%%%%%%%%%%%%%%%%%%%%%%%%%%%%%%%%%%%%%%%%%%%%%%%%%%%%%%%%%%%%%%%%%%%%%%%%%%%%%%%%%%%%%%%%%%%%%%%%%%%%%%%%%%%%%%%%%%%%%%%%%%%%%%%%%%%%%%%%%%%%%%%%%%%%%%%%%%%%%%%%%%%%%%%%%%%%%%%%
%%%%%%%%%%%%%%%%%%%%%%%%%%%%%%%%%%%%%%%%%%%%%%%%%%%%%%%%%%%%%%%%%%%%%%%%%%%%%%%%%%%%%%%%%%%%%%%%%%%%%%%%%%%%%%%%%%%%%%%%%%%%%%%%%%%%%%%%%%%%%%%%%%%%%%%%%%%%%%%%%%%%%%%%%%%%%%%%%%%%%%%%%%%%%%%%%%%%%%%%%%%%%%%%%%%%%%%%
%%%%%%%%%%%%%%%%%%%%%%%%%%%%%%%%%%%%%%%%%%%%%%%%%%%%%%%%%%%%%%%%%%%%%%%%%%%%%%%%%%%%%%%%%%%%%%%%%%%%%%%%%%%%%%%%%%%%%%%%%%%%%%%%%%%%%%%%%%%%%%%%%%%%%%%%%%%%%%%%%%%%%%%%%%%%%%%%%%%%%%%%%%%%%%%%%%%%%%%%%%%%%%%%%%%%%%%%

\mysection{Removable singularities}
Throughout this article we denote by $c$ and $C$ generic constants the value of which may vary from one occurrence to another even within a single string of estimates, and by $c_j$, ($j=1,2,...$) some constants which have a more important significance and a more precise dependence with respect to the parameters. 
\subsection{\bf {\it A priori} estimates}
We give two estimates for positive solutions of $(\ref{Z1})$ which differ according to the sign of $M$. If $G$ is an open subset of $\BBR^N$ we set $d_{_G}(x)=\dist (x,\prt G)$
%%%%%%%%%%%%%%%%%%%%%%%%%%%%%%%%%%%%%%%%%%%%%%%%%%%%%%%%%%%%%%%%%%%%%%%%%%%%%%%%%%%%%%%%%%%%%%%%%%%%%%%%%%%%%%%%%%%%%%%%%%%%%%%%%%%%%%%%%%%%%%%%%%%%%%%%%%%%%%%%%%%%%%%%%%%%%%%%%%%%%%%%%%%%%%%%%%%%%%%%%%%%%%%%%%%%%%%%
\bprop {dom} Let $G\subset\BBR^N$ be an open subset, $M> 0$ and $1<q<p$. If $u\in C^{1}(G)$ is a nonnegative solution of 
$(\ref{Z1})$, there holds,
\begin{equation}\label{1Q1}
u(x)\leq c_5\max\left\{M^{\frac{1}{p-q}}(d_{_G}(x))^{-\frac{q}{p-q}},(d_{_G}(x))^{-\frac{2}{p-1}}\right\}\quad\text{for all }\,x\in G,
\end{equation}
for some $c_5=c_5(N,p,q)>0$.
\es
%%%%%%%%%%%%%%%%%%%%%%%%%%%%%%%%%%%%%%%%%%%%%%%%%%%%%%%%%%%%%%%%%%%%%%%%%%%%%%%%%%%%%%%%%%%%%%%%%%%%%%%%%%%%%%%%%%%%%%%%%%%%%%%%%%%%%%%%%%%%%%%%%%%%%%%%%%%%%%%%%%%%%%%%%%%%%%%%%%%%%%%%%%%%%%%%%%%%%%%%%%%%%%%%%%%%%%%%
\Proof 
Following the method of Keller \cite{Ke} and Osserman \cite{Os}, we fix $x\in G$ and $0<a<d_{_G}(x)$ , and introduce $U(z)=\gl(a^2-|z-x|^2)^{-b}$ for some 
$b>0$. Then putting $r=|x-z|$ and $\tilde U(r)=U(z)$, we have in $B_a(x)$
$$\BA{lll}
L\tilde U=-\tilde U''-\myfrac{N-1}{r}\tilde U'-M|\tilde U'|^q+\tilde U^p\\[4mm]
\phantom{L\tilde U}
=\gl(a^2-r^2)^{-2-b}\left[\gl^{p-1}(a^2-r^2)^{2-b(p-1)}+2b(N-2(b+1))r^2-2Nba^2\right.\\[4mm]
\phantom{L\tilde U-------------}
\left.-M2^qb^{q}\gl^{q-1}r^q(a^2-r^2)^{2+b-q(b+1)}\right].
\EA
$$
If $M>0$, the two necessary conditions on $b>0$ to be fulfilled is order that $\tilde U$ is a supersolution in $B_{|a|}(a)$ are 
$$\BA {lll}&(i)\qquad\qquad 2-b(p-1)\leq 0\Longleftrightarrow b(p-1)\geq 2,\qquad\qquad\qquad\qquad\\[4mm]
&(ii)\qquad\qquad 2+b-q(b+1)\geq 2-b(p-1)\Longleftrightarrow b(p-q)\geq q.\qquad\qquad\qquad\qquad
\EA$$
The above inequalities are satisfied if
\begin{equation}\label{1Q3a}b= \max\left\{\myfrac{2}{p-1},\myfrac{q}{p-q}\right\}.
\end{equation}
If $q>\frac{2p}{p+1}$ then $b=\frac{q}{p-q}$ and 
$$L\tilde U\geq  \gl\left(a^2-r^2\right)^{-\frac{2p-q}{p-q}}\left[\gl^{q-1}\left(\gl^{p-q}-M2^qb^{q}\gr^q\right)\left(a^2-r^2\right)^{\frac{2p-q(p+1)}{p-q}}-(3b+1)Na^2\right].
$$
There exists $c^1_5>0$ depending on $N$, $p$ and $q$ such that if we choose
$$\gl= c^1_5\max\left\{M^{\frac{1}{p-q}}a^{\frac{q}{p-q}},a^{\frac{2p(q-1)}{(p-1)(p-q)}}\right\},
$$
there holds
\begin{equation}\label{1Q4}
L\tilde U\geq 0\quad\text{in }B_a(x).
\end{equation}
Since $\tilde U(z)\to\infty$ when $r\to a$, we derive by the maximum principle that $u\leq \tilde U$ in $B_a(x)$. In particular 
\begin{equation}\label{1Q5}
u(x)\leq \tilde U(x)=\gl a^{-\frac{2q}{p-q}}=c^1_5\max\left\{M^{\frac{1}{p-q}}a^{-\frac{q}{p-q}},a^{-\frac{2}{p-1}}\right\}.
\end{equation}
If $q\leq\frac{2p}{p+1}$ then $b=\frac{2}{p-1}$ and
$$\BA{lll}L\tilde U\geq  \gl\left(|a|^2-\gr^2\right)^{-\frac{2p}{p-1}}\left[\gl^{p-1}+\myfrac{2}{p-1}\left(N-\myfrac{2(p+1)}{p-1}\right)\gr^2-\myfrac{2N}{p-1}|a|^2\right.\\[4mm]\phantom{------------------}
\left.-M2^q\left(\myfrac{2}{p-1}\right)^q\gl^{q-1}\gr^q\left(|a|^2-\gr^2\right)^{\frac{2p-q(p+1)}{p-1}}\right]\\[4mm]\phantom{L\tilde U}
\geq  \gl\left(|a|^2-\gr^2\right)^{-\frac{2p}{p-1}}\left[\gl^{p-1}-c_2|a|^2-c_3\gl^{q-1}M|a|^{\frac{4p-q(p+3)}{p-1}}\right].
\EA$$
Hence, if $q=\frac{2p}{p+1}$, $(\ref{1Q4})$ holds if for some $c^2_5>0$ depending on $N,p,q$, 
$$\gl= c^2_5\max\left\{M^{\frac{p+1}{p(p-1)}},1\right\}|a|^{\frac{2}{p-1}},
$$
which yields
\begin{equation}\label{1Q6}
u(x)\leq \tilde U(x)=\gl a^{-\frac{4}{p-1}}=c^2_5\max\left\{M^{\frac{p+1}{p(p-1)}},1\right\}a^{-\frac{2}{p-1}}.
\end{equation}
While if $q<\frac{2p}{p+1}$, we choose 
$$\gl= c^3_5\max\left\{M^{\frac{1}{p-q}}a^{\frac{4p-q(p+3)}{(p-1)(p-q)}},a^{\frac{2}{p-1}}\right\},
$$
where $c^3_5>0=c^3_5(N,p,q)$, which implies
\begin{equation}\label{1Q6a}
u(x)\leq \tilde U(x)=\gl a^{-\frac{4}{p-1}}=c^3_5\max\left\{M^{\frac{1}{p-q}}a^{-\frac{q}{p-q}},a^{-\frac{2}{p-1}}\right\}.
\end{equation}
By letting $a\uparrow d_{_G}(x)$ we derive $(\ref{1Q1})$ with a constant $c_5=c_5^3$, depending on $N,p,q$.\smallskip
\qeda
%%%%%%%%%%%%%%%%%%%%%%%%%%%%%%%%%%%%%%%%%%%%%%%%%%%%%%%%%%%%%%%%%%%%%%%%%%%%%%%%%%%%%%%%%%%%%%%%%%%%%%%%%%%%%%%%%%%%%%%%%%%%%%%%%%%%%%%%%%%%%%%%%%%%%%%%%%%%%%%%%%%%%%%%%%%%%%%%%%%%%%%%%%%%%%%%%%%%%%%%%%%%%%%%%%%%%%%%
\bcor{dom+} Under the assumptions of \rprop{dom} with $G=B_{2R}\setminus\{0\}$, there holds for $x\in B_{R}\setminus\{0\}$,\smallskip
\begin{equation}\label{1Q12}\BA {lll}
 u(x)\leq c_5\max\left\{M^{\frac{1}{p-q}}|x|^{-\frac{q}{p-q}},|x|^{-\frac{2}{p-1}}\right\}.
\EA\end{equation}
\es

%%%%%%%%%%%%%%%%%%%%%%%%%%%%%%%%%%%%%%%%%%%%%%%%%%%%%%%%%%%%%%%%%%%%%%%%%%%%%%%%%%%%%%%%%%%%%%%%%%%%%%%%%%%%%%%%%%%%%%%%%%%%%%%%%%%%%%%%%%%%%%%%%%%%

We infer from \rprop{dom} an estimate of the gradient of a positive solution when $M>0$. We set $\gs=2p-q(p+1)$, then $\gs>0$ (resp. $\gs<0$) according 
$2p>q(p+1)$ (resp. $2p<q(p+1)$). 

\bprop{domgrad} Let $p>q>1$. For any $M_0>0$ and $R>0$ there exists a constant $c_8=c_8(N,p,q,M_0R^{\frac{\gs}{p-1}})$ such that, for $0<M\leq M_0$ there holds:\smallskip

\nind (i) If  $q\leq \frac{2p}{p+1}$ (then $\gs\geq 0$), any  positive solution  $u$ of $(\ref{Z1})$ in $B_{2R}\setminus\{0\}$ satisfies
\begin{equation}\label{1Q1g}\displaystyle
|\nabla u(x)|\leq c_8\max\left\{M^{\frac{1}{p-q}}|x|^{-\frac{p}{p-q}},|x|^{-\frac{p+1}{p-1}}\right\},
\end{equation}
for all $x\in B_{R}\setminus\{0\}$.
\smallskip

\nind (ii) If $\frac{2p}{p+1}\leq q\leq 2$ (then $\gs\leq 0$), any  positive solution  $u$ of $(\ref{Z1})$ in $\BBR^N\setminus\overline B_{\frac R2}$ satisfies $(\ref{1Q1})$ 
for all $x\in \BBR^N\setminus B_{R}$.

\es
\Proof (i) For $0<r<2R$ we set 
$$u(x)=r^{-\frac{2}{p-1}}u_r(\tfrac{x}{r})=r^{-\frac{2}{p-1}}u_r(y)\quad\text{with }\,y=\tfrac{x}{r}.
$$
If $\frac{r}{2}<|x|<2r$, then $\frac{1}{2}<|y|<2$ and $u_r>0$ satisfies
\begin{equation}\label{1Q15m}-\Gd u_r+u_r^p-Mr^{\frac{2p-q(p+1)}{p-1}}|\nabla u_r|^q=0\qquad\text{in }\; B_2\setminus B_{\frac 12}.
\end{equation}
Since $0<Mr^{\frac{\gs}{p-1}}\leq M(2R)^{\frac{\gs}{p-1}}\leq M_0(2R)^{\frac{\gs}{p-1}}$ as $\gs\geq 0$, it follows that 
\begin{equation}\label{1Q15n}\BA {lll}
\max\left\{|\nabla u_r(z)|: \frac{2}{3}<|z|<\frac{3}{2}\right\}\leq c\max\left\{|u_r(z)|: \frac{1}{2}<|z|<2\right\},
\EA\end{equation}
where $c$ depends on $N,p,q$ and $R^{\frac{\gs}{p-1}}M_0$ (see e.g. \cite[Chapter 13]{GT}). From \rprop{dom} there holds
$$
\max\left\{|u_r(z)|: \tfrac{1}{2}<|z|\leq 2\right\}\leq 2^{\frac{2}{p-1}}c_5\max\left\{M^{\frac{1}{p-q}}r^{\frac{2p-q(p+1)}{(p-1)(p-q)}},1\right\}
$$
 by $(\ref{1Q1})$. Therefore
\begin{equation}\label{1Q16-}\BA {lll}
\max\left\{|\nabla u(y)|: \tfrac{r}{2}<|z|<2r\right\}\leq 2^{\frac{2}{p-1}}cc_5r^{-\frac{p+1}{p-1}}\max\left\{M^{\frac{1}{p-q}}r^{\frac{2p-q(p+1)}{(p-1)(p-q)}},1\right\}
\\[2mm]\phantom{\max\left\{|\nabla u(y)|: \tfrac{r}{2}<|z|<2r\right\}}
\leq c_8\max\left\{M^{\frac{1}{p-q}}|x|^{-\frac{p}{p-q}},|x|^{-\frac{p+1}{p-1}}\right\},
\EA\end{equation}
which implies $(\ref{1Q1g})$. \smallskip 

\nind (ii) For $r>R$ we define $u_r$ as in (i). It satisfies $(\ref{1Q15m})$ and since $\gs\leq 0$, we have again $0<Mr^{\frac{\gs}{p-1}}\leq MR^{\frac{\gs}{p-1}}\leq M_0R^{\frac{\gs}{p-1}}$ if $r\geq R$. Since $1<q<2$, $(\ref{1Q15n})$ holds and we derive $(\ref{1Q1g})$.\qeda\medskip

\nind\Remark If $q=\frac{2p}{p+1}$ the constant $c_8$ depends only on $N$ and $p$.\medskip
 
 The previous estimate necessitates $1<q\leq 2$. This limitation can be by passed in some cases using the Bernstein  approach.

\blemma{domgradlem} Assume $p,q>1$ and $M>0$. If $u\in C^2(\overline B_{2R})$ is a nonnegative solution of  $(\ref{Z1})$ in $B_{2R}$, there holds
\begin{equation}\label{1Q14}\displaystyle
|\nabla u(x)|\leq c_9\left(|x|^{-\frac{1}{q-1}}+\max_{|z-x]\leq\frac{|x|}{2}} u^\frac{p}{q}(z)\right)\qquad\text{for all }\,x\in B_{\frac R2},
\end{equation}
$c_{9}>0$ depends on $N$, $p$, $q$ and $M$.
\es
\Proof Set $z=|\nabla u|^2$, then by a classical computation and the use of Schwarz inequality,
$$-\Gd|\nabla u|^2+\myfrac{1}{N}(\Gd u)^2+\langle\nabla \Gd u,\nabla u\rangle\leq 0.
$$
Replacing $\Gd u$ by its expression from $\CL^{_{^M}}_{p,q}u=0$, we obtain
$$-\Gd z+\myfrac{2}{N}\left(u^{2p}+M^2z^q-2Mu^pz^{\frac q2}\right)+2pu^{p-1}z\leq qMz^{\frac q2-1}\langle\nabla z,\nabla u\rangle.
$$
We notice that
$$\BA {lll}
qMz^{\frac q2-1}\langle\nabla z,\nabla u\rangle\leq qMz^{\frac q2-1}|\nabla z|\sqrt z=qMz^{\frac q2}\myfrac{|\nabla z|}{\sqrt z}\leq \myfrac{M^2z^q}{2N}+\myfrac{2Nq^2}{M^2}\myfrac{|\nabla z|^2}{z},
\EA$$
and
$$\myfrac{4M}{N}u^pz^{\frac q2}\leq \myfrac{M^2z^q}{2N}+\myfrac{8u^{2p}}{NM^2},
$$
thus
$$-\Gd z+\myfrac{M^2z^q}{N}\leq \myfrac{2Nq^2}{M^2}\myfrac{|\nabla z|^2}{z}+\myfrac{2}{N}\left(\myfrac{4}{M^2}-1\right)u^{2p}.
$$
For simplicity we set 
$$\displaystyle A=\myfrac{M^2}{N}\,,\; B=\myfrac{2Nq^2}{M^2}\,\text{and }\, C=\myfrac{2}{N}\left(\myfrac{4}{M^2}-1\right)_+\max_{|z-x|\leq \frac{|x|}{2}}u^{2p}(z)
$$
Then $z$ satisfies 
$$-\Gd z+Az^q\leq B\myfrac{|\nabla z|^2}{z}+C\quad\text{in }B_{\frac R2}(x)$$
and obtain by \cite[Lemma 3.1]{BV} (see also a simpler approach in \cite[Lemma 2.2]{BVGHV2}),
\begin{equation}\label{1Q15}\displaystyle
z(x)\leq c_{10}\left(|x|^{-\frac{2}{q-1}}+C^{\frac 1q}\right)
\end{equation}
where $c_{10}>0$ depends on $N$, $p$, $q$ and $M$. This yields $(\ref{1Q14})$.\qeda\medskip

\nind\Remark The constants $c_{9}$ and $c_{10}$ can be expressed in terms of $M$, but their stability when $M\to 0$ is not clear since in the limit case of 
the equation $\CL_pu=0$ the estimate of the gradient obtained by a very different and much simpler method combining the Keller-Osserman estimate and scaling methods.  \medskip

Using \rcor{dom+} we obtain the new estimate
\bcor{domgradlem2} Assume $1<q<p$ and $M>0$. Then any nonnegative  solution $u\in C^2(B_{2R})$ of $(\ref{Z1})$ satisfies
\begin{equation}\label{1Q16}\displaystyle
|\nabla u(x)|\leq c_{11}\left(|x|^{-\frac{1}{q-1}}+\max\left\{M^{\frac{p}{q(p-q}}|x|^{-\frac{p}{p-q}}, |x|^{-\frac{2p}{q(p-1)}}\right\}\right)\quad\text{for all }\,x\in B_{\frac R2},
\end{equation}
where $c_{11}>0$ depends on $N$, $p$, $q$ and $M$.
\es

Then we can combine this estimate with  \rprop{dom} to complete the cases not treated in \rprop{domgrad}.

\bprop{domgrad2} Let $1<q<p$. For any $M>0$ there exists a constant $c_{12}=c_{12}(N,p,q,M)>0$ such that:\smallskip

\nind (i) If  $\frac{2p}{p+1}\leq q<p$, any  positive solution  $u$ of $(\ref{Z1})$ in $B_{2R}\setminus\{0\}$ with $0<R\leq 1$ satisfies,
\begin{equation}\label{1Q17}\displaystyle
|\nabla u(x)|\leq c_{12}\max\left\{M^{\frac{p}{q(p-q)}}|x|^{-\frac{p}{p-q}},|x|^{-\frac{2p}{q(p-1)}}\right\}\quad \text{for all } x\in B_{R}\setminus\{0\}.
\end{equation}

\nind (ii) If $1<q\leq \frac{2p}{p+1}$, any  positive solution  $u$ of $(\ref{Z1})$ in $\BBR^N\setminus\overline B_{\frac R2}$ with $R\geq 1$ satisfies,
\begin{equation}\label{1Q18}\displaystyle
|\nabla u(x)|\leq c_{12}\max\left\{M^{\frac{p}{q(p-q)}}|x|^{-\frac{p}{p-q}},|x|^{-\frac{1}{q-1}}\right\}\quad \text{for all } x\in \BBR^N\setminus B_{R}.
\end{equation}
\es
\Proof We can compare the different exponents of $|x|$ which appear in the expressions $(\ref{1Q1g})$ and $(\ref{1Q16})$
\begin{equation}\label{1Q19}\BA {lll}\displaystyle
(i)\quad&\myfrac{p}{p-q}<\myfrac{p+1}{p-1}<\myfrac{2p}{q(p-1)}<\myfrac{1}{q-1}\quad&\text{if }1<q<\myfrac{2p}{p+1},\quad\\[4mm]
(ii)\quad&\myfrac{1}{q-1}<\myfrac{2p}{q(p-1)}<\myfrac{p+1}{p-1}<\myfrac{p}{p-q}&\text{if }\myfrac{2p}{p+1}<q<p,\quad
\EA\end{equation}
with equality if $q=\frac {2p}{p+1}$. If $1<q\leq \frac{2p}{p+1}$ (resp. $\frac{2p}{p+1}\leq \leq 2$), estimate  $(\ref{1Q1g})$  is better than $(\ref{1Q16})$ in $B_R\setminus\{0\}$ (resp. $\BBR^N\setminus B_{2R}$). Then $(\ref{1Q17})$ and $(\ref{1Q18})$ follow from $(\ref{1Q16})$ and $(\ref{1Q19})$.\qeda
\medskip

In the case $M< 0$ an upper estimate on a solution is obtained by combining a result of Lions and the method of Keller and Osserman. 
\bprop {dom2} Let $G\subset\BBR^N$ be an open subset, $M\leq 0$ and $p,q>1$. If $u\in C^{1}(G)$ is a nonnegative solution of 
$\CL^{_{^M}}_{p,q} u=0$, there exists $c_6=c_6(N,p)>0$, $c_7=c_7(N,q)>0$ and $\gd=\gd(G)>0$ such that there holds for all $x\in G$ and 
$\O<\gd\leq \gd(G)$,
\begin{equation}\label{1Q2}\displaystyle
u(x)\leq \min\left\{c_6(d_{_G}(x))^{-\frac{2}{p-1}}, c_7|M|^{-\frac{1}{q-1}}(d_{_G}(x))^{-\frac{2-q}{q-1}}+\max_{d_{_G}(z)=\gd} u(z),\right\}.
\end{equation}
\es
\Proof This estimates follows from the fact that the solutions of $\CL^{_{^M}}_{p,q} u=0$ are subsolutions of $\CL_{p} u=0$ and $\CR^{_{^M}}_{q} u=0$. The estimate $u(x)\leq c_6(d_{_G}(x))^{-\frac{2}{p-1}}$ corresponds to the Keller-Osserman estimate for solutions of $\CL_{p} u=0$. The second estimate 
corresponds to the fact that if $u$ is a positive solution of  $\CR^{_{^M}}_{q} u=0$ in $G$ there holds (see \cite[Theorem IV-1]{Li})
$$|\nabla u(x)|\leq c'_7|M|^{-\frac 1{q-1}}(d_{_G}(x))^{-\frac{1}{q-1}}.
$$
Integrating this inequality yields the second part of the inequality. \qeda\medskip

\nind\Remark This estimate can be transformed into the universal estimate
\begin{equation}\label{1Q2'}\displaystyle
u(x)\leq \min\left\{c_6(d_{_G}(x))^{-\frac{2}{p-1}}, c_7|M|^{-\frac{1}{q-1}}(d_{_G}(x))^{-\frac{2-q}{q-1}}+c_6\gd^{-\frac{2}{p-1}},\right\},
\end{equation}
since $\displaystyle\max_{d_{_G}(z)=\gd} u(z)\leq c_6\gd^{-\frac{2}{p-1}}$ by $(\ref{1Q2})$.\medskip

The gradient estimates  are due to Nguyen \cite[Proposition 1.1]{NPT}. Below we recall his result proved by the Bernstein method in a more general framework  but which can also be obtained by scaling techniques in the present case.
\bprop{domgrad3} Let $p>1$ and $1<q<2$. For any $M<0$ and $R>0$ there exists a constant $c'_{12}=c'_{12}(N,p,q,M,)>0$ such that: if $u$ is a
 positive solution of $(\ref{Z1})$ in $B_{2R}\setminus\{0\}$, there holds
 \begin{equation}\label{1Q20}\displaystyle
u(x)+|x||\nabla u(x)|\leq c'_{12}\max\left\{|x|^{-\frac{2}{p-1}},|x|^{-\frac{2-q}{q-1}}\right\}\quad \text{for all } x\in B_{R}\setminus\{0\}.
\end{equation}
\es

\nind \Remark There are many estimates of positive solutions of $(\ref{Z1})$ (or even with $u^p$ replaced by $f(u)$) in a domain which tends to infinity on the boundary ({\it large solutions}) or of solutions in $\BBR^N$ ({\it ground states}). Many of  these estimates are obtained by comparison with one dimensional problems and they can be found in  \cite{AGQ}, \cite{BG}, \cite{FQS}.
%%%%%%%%%%%%%%%%%%%%%%%%%%%%%%%%%%%%%%%%%%%%%%%%%%%%%%%%%%%%%%%%%%%%%%%%%%%%%%%%%%%%%%%%%%%%%%%%%%%%%%%%%%%%%%%%%%%%%%%%%%%%%%%%%%%%%%%%%%%%%%%%
%%%%%%%%%%%%%%%%%%%%%%%%%%%%%%%%%%%%%%%%%%%%%%%%%%%%%%%%%%%%%%%%%%%%%%%%%%%%%%%%%%%%%%%%%%%%%%%%%%%%%%%%%%%%%%%%%%%%%%%%%%%%%%%%%%%%%%%%%%%%%%%%
%%%%%%%%%%%%%%%%%%%%%%%%%%%%%%%%%%%%%%%%%%%%%%%%%%%%%%%%%%%%%%%%%%%%%%%%%%%%%%%%%%%%%%%%%%%%%%%%%%%%%%%%%%%%%%%%%%%%%%%%%%%%%%%%%%%%%%%%%%%%%%%%
%%%%%%%%%%%%%%%%%%%%%%%%%%%%%%%%%%%%%%%%%%%%%%%%%%%%%%%%%%%%%%%%%%%%%%%%%%%%%%%%%%%%%%%%%%%%%%%%%%%%%%%%%%%%%%%%%%%%%%%%%%%%%%%%%%%%%%%%%%%%%%%%

\subsection{Proof of \rth{Remov-1}}
Without loss of generality we can assume that $u\in C^2(\overline\Gw\setminus\{0\}$ and $\overline B_{2R_0}\subset\Gw$ with $2R_0\leq 1$. If $M\leq 0$, $u$ is a nonnegative subsolution of 
$-\Gd u+u^p=0$, hence it is bounded in $\overline\Gw$ by \cite{BrVe}. \smallskip

\noindent {\it Step 1. }We assume $M>0$ and we prove first that under condition (i) or (ii), $|\nabla u|^q\in L^1(\Gw)$,  $u\in L^p(\Gw)$, and then
   \begin{equation}\label{RW231}\BA {lll}
\myint{\Gw}{}\left(-u\Gd\gz+u^p\gz-M|\nabla u|^q\gz\right)dx=0\qquad\forall \gz\in W^{2,\infty}(\Gw)\cap C^1_c(\overline\Gw).
\EA\end{equation}
By \rprop{domgrad} 
$$|\nabla u|^q\leq c|x|^{-\frac{(p+1)q}{p-1}}\quad\text{in }B_{R_0},
$$
since $q\leq \frac{2p}{p+1}$, and where $c$ depends also on $M$. By (i) or (ii), $\frac{(p+1)q}{p-1}<N$. Hence $\nabla u\in L_{loc}^q(\Gw)$.
\medskip

For any $\ge>0$ small enough we denote by $\gr_\ge$ a nonnegative  $C^\infty$-function such that  supp($\gr_\ge$)$\,\subset \overline B_\ge$, $0\leq\gr_\ge\leq 1$, $|\nabla\gr_\ge|\leq 2\ge^{-1}\chi_{_{\overline B_\ge}}$ and we set $\eta_\ge=1-\gr_\ge$. Then
  \begin{equation}\label{RW131}\BA {lll}-\myint{B_{2R_0}}{}\langle\nabla u,\nabla\gr_\ge\rangle dx+\myint{B_{2R_0}}{}u^p\eta_\ge dx+
  \myint{\prt B_{2R_0}}{}\myfrac{\prt u}{\prt{\bf n}}dS=M\myint{B_{2R_0}}{}|\nabla u|^q\eta_\ge dx.
\EA\end{equation}
Next
  \begin{equation}\label{RW132}\left|\myint{B_{2R_0}}{}\langle\nabla u,\nabla\gr_\ge\rangle  dx\right|\leq 2c_N\ge^{\frac{N}{q'}-1}\left(\myint{B_{\ge}}{}|\nabla u|^q dx\right)^{\frac1q}\to 0\quad\text{as }\,\ge \to 0,
\end{equation}
since $1<q\leq \frac{N}{N-1}$. Since $|\nabla u|^q\in L^1(B_{2R_0})$ we deduce by monotone convergence that $u^p\in L^1(B_{2R_0})$. Finally, if  $\gz\in C^\infty_0(\Gw)$ and $\gz_\ge=\gz\eta_\ge$, we have
$$\myint{\Gw}{}\langle\nabla u,\nabla\gz_\ge\rangle dx+\myint{\Gw}{}u^p\gz_\ge dx
-M\myint{\Gw}{}|\nabla u|^q\gz_\ge dx=0.
$$
Letting $\ge\to 0$ and using $(\ref{RW132})$, we infer that $u$ satisfies
$$\myint{\Gw}{}\langle\nabla u,\nabla\gz\rangle dx+\myint{\Gw}{}u^p\gz dx
-M\myint{\Gw}{}|\nabla u|^q\gz dx=0.
$$
Hence it is a weak solution of $(\ref{Z11})$ in $\Gw$. 
\smallskip

\noindent {\it Step 2. } Let us assume that $p>\frac{N}{N-2}$. If $u$ is nonnegative and not identically zero, it is positive in $\Gw\setminus\{0\}$ by the maximum principle. We set $u=v^b$ with $0<b\leq 1$. Then 
 \begin{equation}\label{RW1}
\displaystyle
-\Gd v-(b-1)\myfrac{|\nabla v|^2}{v}+\myfrac{1}{b}v^{1+(p-1)b}-Mb^{q-1}v^{(b-1)(q-1)}|\nabla v|^q=0.
\end{equation}
 For $\ge>0$, 
 $$v^{(b-1)(q-1)}|\nabla v|^q\leq \myfrac{q\ge^{\frac{2}{q}}}{2}\myfrac{|\nabla v|^2}{v}+\myfrac{2-q}{2\ge^{\frac{2}{2-q}}}v^{1+\frac{2b(q-1)}{2-q}}.
 $$
 Therefore
  \begin{equation}\label{RW2}
\displaystyle
-\Gd v+\left(1-b-M\myfrac{qb^{q-1}\ge^{\frac{2}{q}}}{2}\right)\myfrac{|\nabla v|^2}{v}+\myfrac{1}{b}v^{1+b(p-1)}-Mb^{q-1}\myfrac{2-q}{2\ge^{\frac{2}{2-q}}}v^{1+\frac{2b(q-1)}{2-q}}=0.
\end{equation}
We notice that $1+\frac{2b(q-1)}{2-q}=1+b(p-1)-a$ with $a=b\frac{2p-(p+1)q}{2-q}\geq 0$. We fix $b$ as follows,
  \begin{equation}\label{RW3}(p-1)b+1=\frac{N}{N-2}\Longleftrightarrow b=\frac{2}{(N-2)(p-1)},
\end{equation}
 hence $p> \frac{N}{N-2}$ if and only if $0<b< 1$. Next we impose
  \begin{equation}\label{RW4}
1-b-M\myfrac{qb^{q-1}\ge^{\frac{2}{q}}}{2}=0\Longleftrightarrow\ge=\left(\frac{2(1-b)}{Mqb^{q-1}}\right)^{\frac q2}=\left(\frac{2((N-2)p-N)}{Mqb^{q-1}(N-1)(p-1)}\right)^{\frac q2}\!\!\!.
\end{equation}
This transforms $(\ref{RW2})$ into
  \begin{equation}\label{RW5}
\displaystyle
-\Gd v+\myfrac{(N-2)(p-1)}{2}v^{\frac{N}{N-2}}-\myfrac{(2-q)b^{q-1}}{2}\left(\myfrac{q}{2(1-b)}\right)^{\frac{q}{2-q}}M^{\frac{2}{2-q}}v^{\frac{N}{N-2}-a}\leq 0.
\end{equation}

\noindent We first assume that $0<q<\frac{2p}{p+1}$. Then  $a>0$, hence
there exists $A>0$, depending on $M$, such that 
  \begin{equation}\label{RW6}
\displaystyle
-\Gd v+\myfrac{(N-2)(p-1)}{4}v^{\frac{N}{N-2}}\leq A.
\end{equation}
Set $\tilde v=(v-cA^{\frac{N-2}{N}})_+^{\frac{N}{N-2}}$ with $c=\left(\frac{4}{(N-2)(p-1)}\right)^{\frac{N}{N-2}}$ satisfies
  \begin{equation}\label{RW7}
\displaystyle
-\Gd \tilde v+\myfrac{(N-2)(p-1)}{4}\tilde v^{\frac{N}{N-2}}\leq 0.
\end{equation}
By \cite{BrVe}, $\displaystyle\tilde v\leq \max_{\prt\Gw} \tilde v$ which implies $\displaystyle v\leq cA^{\frac{N-2}{N}}+\max_{\prt\Gw} v$ and therefore $u(x)\leq B$ for some $B\geq 0$ in $\Gw$.  Furthermore $|\nabla u|^{q-1}\in L^{\frac{q}{q-1}}(\Gw)$ since $\nabla u\in L^q(\Gw)$, and $\frac{q}{q-1}>N$ as we assume $q<\frac{N}{N-1}$.
Writing $(\ref{Z11})$ under the form
$$-\Gd u+u^p-MC(x)|\nabla u|=0,
$$
with $C(x)=|\nabla u(x)|^{q-1}$, it follows from Serrin's theorem \cite[Theorem 10]{SerActa} that the singularity at $0$ is removable and $u$ can be extended as a 
regular solution of $(\ref{Z11})$ in $\Gw$. Hence $ u\in C^2(\Gw)$. \smallskip

\noindent Then we assume that $q=\frac{2p}{p+1}$. By the choice of $b$ in  $(\ref{RW3})$, inequality  $(\ref{RW2})$ becomes 
  \begin{equation}\label{RW11}
\displaystyle
-\Gd v+\left(1-b-\myfrac{Mpb^{\frac{p-1}{p+1}}\ge^{\frac{p+1}{p}}}{p+1}\right)\myfrac{|\nabla v|^2}{v}+\left(\myfrac{1}{b}-
\myfrac{Mb^{\frac{p-1}{p+1}}}{(p+1)\ge^{p+1}}\right)v^{\frac{N}{N-2}}\leq 0.
\end{equation}
Notice that 
  \begin{equation}\label{RW11*}\myfrac{1}{b}-
\myfrac{Mb^{\frac{p-1}{p+1}}}{(p+1)\ge^{p+1}}=0\Longleftrightarrow \ge=\left(\frac{M}{p+1}\right)^{\frac{1}{p+1}}
b^{\frac{2p}{(p+1)^2}},
\end{equation}
and therefore
  \begin{equation}\label{RW11**}1-b-\myfrac{Mpb^{\frac{p-1}{p+1}}\ge^{\frac{p+1}{p}}}{p+1}=1-b-pb\left(\frac{M}{p+1}\right)^{\frac{p+1}p}.
\end{equation}
This coefficient vanishes if 
$$p\left(\frac{M}{p+1}\right)^{\frac{p+1}p}=\myfrac{p(N-1)-(N+1)}{2}.
$$
Therefore, if $M$ satisfies 
  \begin{equation}\label{RW12}
\displaystyle
p\left(\frac{M}{p+1}\right)^{\frac{p+1}p}=\myfrac{p(N-2)-N}{2},
\end{equation}
we can choose $\ge>0$ so that the coefficient of $v^{(p-1)b+1}$ in $(\ref{RW12})$ is equal to some $\gt>0$. Therefore $v$ satisfies 
  \begin{equation}\label{RW13}\BA {lll}
-\Gd v+\gt v^{\frac{N}{N-2}}\leq 0\qquad\text{in }\,\Gw\setminus\{0\}.
\EA\end{equation}
It follows by \cite{BrVe}, $v\leq \displaystyle\max_{\prt\Gw}v$ and the same type of uniform estimate holds for $u$. This ends the case $p>\frac{N}{N-2}$.\smallskip

%%%%%%%%%%%%%%%%%%%%%%%%%%%%%%%%%%%%%%%%%%%%%%%%%%%%%%%%%%%%%%%%%%%%%%%%%%%%%%%%%%%%%%%%%%%%%%%%%%%%%%%%%%%%%%%%%%%%%%%%%%%%%%%%%%%%%%%%%%%%%%%%

\noindent {\it Step 3. } Finally we assume $p=\frac{N}{N-2}$ and $1<q<\frac{2p}{p+1}=\frac{N}{N-1}$. From $(\ref{1Q14})$,  
$$M|\nabla u(x)|^q\leq c_{9}|x|^{-q\frac{p+1}{p-1}}=c_{9}|x|^{-q(N-1)}:=Q(x).$$
and $Q\in L^1(B_{2R_0})$. Let $\{\gs_n\}\subset C^\infty_0(\BBR^N)$ such that $0\leq \gs_n\leq 1$
$$\gs_n(x)=\left\{\BA {lll}1\qquad&\text{if }\frac{2}{n}\leq |x|\leq R_0\\[2mm]
0\qquad&\text{if }|x|\in [0, \frac{1}{n}]\cup[2R_0,\infty), 
\EA\right.
$$
and 
$$|\Gd\gs_n|\leq 2Nn^2\chi_{_{B_{\frac 2n}\setminus B_{\frac 1n}}}+\phi,$$ 
where $\phi$ is a smooth nonnegative function with support in $B_{2R_0}\setminus B_{R_0}$. Then
  \begin{equation}\label{RW14*}-\myint{\{\frac 1n\leq|x|\leq \frac 2n\}}{}\!\!\!\!\!\!\!\!\!\!\!\!\!\!\!\!u\Gd\gs_ndx-\myint{\{R_0\leq|x|\leq 2R_0\}}{}\!\!\!\!\!\!\!\!\!\!\!\!\!\!\!\!u\Gd\gs_ndx+\myint{\frac 1n\leq |x|}{}u^p\gs_ndx=
M\myint{\frac 1n\leq |x|}{}|\nabla u|^q\gs_ndx.
\end{equation}
The right-hand side of $(\ref{RW14*})$ is bounded since $|\nabla u|\in L^q(B_{2R_0})$, the second term on the left is also uniformly bounded. Using the fact that 
$|x|^{N-2}u(x)$ is bounded by  $(\ref{1Q1})$, we get
$$\left|\myint{\{\frac 1n\leq|x|\leq \frac 2n\}}{}\!\!\!\!\!\!\!\!\!\!\!\!u\Gd\gs_ndx\right|\leq C,
$$
for some $C>0$ independent of $n$. Letting $n\to\infty$ we infer that $u\in L^p_{loc}(\Gw)$. By the maximum principle 
  \begin{equation}\label{RW14**}\displaystyle u(x)\leq u_1(x)=C{\bf G}^{B_{2R_0}}[Q](x)+\max_{|z|=2R_0}u(z),
\end{equation}
where ${\bf G}^{B_{2R_0}}$ denotes the Green kernel in $B_{2R_0}$. Since $Q(x)=C|x|^{-q(N-1)}$, a direct computation shows that 
$u_1(x)\leq c_NC|x|^{2-q(N-1)}= c_NC|x|^{2-N+\ge}$ for some $\ge>0$. We can write $(\ref{Z1})$ under the form
$$-\Gd u+c(x)u+d(x)|\nabla u|=0\quad\text{in }\Gw\setminus\{0\},
$$
with $c(x)=u^{\frac {2}{N-2}}$ and $d(x)=|\nabla u|^{q-1}$. Then $c\in L^{\frac{N}{2}+\ge_1}(B_{2R_0})$ and $d\in L^{N+\ge_2}(B_{2R_0})$; with $\ge_1,\ge_2>0$. 
It follows from \cite[Theorem 10]{SerActa} that $0$ is a removable singularity for $u$ in the sense that it can be extended as a $C^2$ solution in $\Gw$. 
\qeda\medskip

 When the conditions of the theorem are not fulfilled there exist singular solutions.  However these singular solutions may exhibit different types of behaviour according $1<q<\frac{2p}{p+1}$, $\frac{2p}{p+1}<q<2$ and  $q=\frac{2p}{p+1}$. In this case there may exist radial separable solutions of $(\ref{Z11})$ under the form $u_X(r)=X r^{-\frac 2{p-1}}$. Setting $\ga=\frac 2{p-1}$, then $X$ satisfies
   \begin{equation}\label{RW15-c}
\Phi_p(X):=X^{p-1}-M\ga^{\frac{2p}{p+1}}X^{\frac{p-1}{p+1}}+\ga(N-2-\ga)=0
\end{equation}
  The following result is easy to prove by a standard analysis of the function $\Phi_p$.
  \bprop{separa} Let $p>1$ and $M\in\BBR$. \smallskip
  
\nind (i) If $M$ is arbitrary and $1<p<\frac{N}{N-2}$, or $M>0$ and $p=\frac{N}{N-2}$,  
there exists one and only one positive solution $X_1$ to $(\ref{RW15-c})$. \smallskip

\nind (ii) If $p>\frac{N}{N-2}$ and  $M>m^*$, 
there exist two positive solutions $X_1<X_2$ to $(\ref{RW15-c})$.\smallskip

\nind (iii) If $p>\frac{N}{N-2}$ and  $M=m^*$ there exists one positive solution $X_1$ to $(\ref{RW15-c})$.\smallskip

\nind (iv) If $p> \frac{N}{N-2}$ and $0<M<m^*$, or $M\leq 0$ and $p\geq \frac{N}{N-2}$, there exists no positive solution to $(\ref{RW15-c})$.
\es 
 
 When $q\neq\frac{2p}{p+1}$ the existence of singular solutions is much involved. It is developed in the subsequent paper \cite{BVGHV3}.\medskip

 \noindent\Remark It is noticeable that in the case $q=\frac{2p}{p+1}$, $p>\frac{N}{N-2}$ and $M\geq m^*$, the equation exhibits a phenomenon which is characteristic of Lane-Emden type equations
  \begin{equation}\label{RW16}
-\Gd u=u^p\quad\text{in }\;B_1\setminus\{0\}.
\end{equation}
If $u$ is nonnegative then there exists $\ga\geq 0$ such that 
  \begin{equation}\label{RW16a}
-\Gd u=u^p+\ga\gd_0\quad\text{in }\;\CD'(B_1).
\end{equation}
If $1<p<\frac{N}{N-2}$ then $\ga$ can be positive, but if $p\geq\frac{N}{N-2}$, then $\ga=0$. This means that the singularity cannot be seen 
in the sense of distributions, however there truly exist singular solutions, e.g. if $p>\frac{N}{N-2}$,
  \begin{equation}\label{RW17a}
u_s(x)=c_{N,p}|x|^{-\frac{2}{p-1}}.
\end{equation}
Here also for $q=\frac{2p}{p+1}$, $p>\frac{N}{N-2}$, the isolated singularities are not seen in the sense of distributions. 
%%%%%%%%%%%%%%%%%%%%%%%%%%%%%%%%%%%%%%%%%%%%%%%%%%%%%%%%%%%%%%%%%%%%%%%%%%%%%%%%%%%%%%%which ends the proof. \qeda%%%%%%%%%%%%%%%%%%%%%%%%%%%%%%%%%%%%%%%%%%%%%%%%%%%%%%%%%%%%%%%%%%%%%%%%%%%%%%%%%%%%%%%%%%%%%%%%%%%%%%%%%%%%%%%%%%%%%%%%%%%%%%%%%%%%%%%%%%%%%%%%%%%%%%%%%%%%%%%%%%%%%%%%%%%%%%%%%%%%%%%%%%%%%%%%%%%%%%%%%
\subsection{Proof of \rth{Remov-2}}
 As in the proof of \rth{Remov-1}, we distinguish according $1<q<\frac{2p}{p+1}$ or $q=\frac{2p}{p+1}$. Without loss of generality we can suppose that $u>0$. We perform the same change of unknown as in the previous theorem putting $u=v^b$, but now we choose $b$ as follows
  \begin{equation}\label{RW16b}(p-1)b+1=r\Longleftrightarrow b=\frac{r-1}{p-1},
\end{equation}
and we first assume that 
  \begin{equation}\label{RW16aa}1-b-M\frac{qb^{q-1}\ge^{\frac2q}}{2}=0\Longleftrightarrow \ge=\left(\frac{2(1-b)}{Mqb^{q-1}}\right)^
\frac q2=\left(\frac{2(p-r)}{Mq(p-1)b^{q-1}}\right)^
\frac q2.
\end{equation}
Hence $(\ref{RW5})$ becomes  
  \begin{equation}\label{RW17}
\displaystyle
-\Gd v+\myfrac{p-1}{r-1}v^{r}-\myfrac{(2-q)b^{q-1}}{2}\left(\myfrac{q}{2(1-b)}\right)^{\frac{q}{2-q}}M^{\frac{2}{2-q}}v^{\frac{(2r-p-1)q+2(p-r)}{(p-1)(2-q)}}\leq 0.
\end{equation}
Condition $r\geq \frac{(2r-p-1)q+2(p-r)}{(p-1)(2-q)}$ is equivalent to $2p-q(p+1)\leq r(2p-q(p+1))$, since $1<r<p$. \par
\noindent Assuming first that $q<\frac{2p}{p+1}$, we obtain from $(\ref{RW17})$
  \begin{equation}\label{RW18a}
\displaystyle
-\Gd v+\myfrac{p-1}{2(r-1)}v^{r}\leq A.
\end{equation}
for some constant $A\geq 0$. Since $cap^{\BBR^N}_{2,r'}(K)=0$ the function $v$ is bounded from \cite{BaPi}  and 
$\displaystyle v\leq cA^{\frac1r}+\max_{\prt\Gw}v$ for some $c>0$, hence $u$ is also uniformly upper bounded in $\Gw$ by some constant $a$. \\
Next we have to show that $\nabla u\in L^q(\Gw)$. Let $\{\gr_n\}$ be a sequence of $C^\infty_0(\Gw)$ nonnegative functions such that 
$0\leq\gr_n\leq 1$, $\gr_n=1$ in a  small enough neighborhood of $K$ and $\norm{\gr_n}_{W^{2,r'}}\to 0$ when $n\to\infty$, and set $\eta_n=1-\gr_n$. Since
$$\myint{\Gw}{}u\Gd\gr_n dx-\myint{\prt\Gw}{}\myfrac{\prt u}{\prt{\bf n}}dS+\myint{\Gw}{}u^p\eta_n dx=M\myint{\Gw}{}|\nabla u|^q\eta_n dx,
$$
and 
$$\left|\myint{\Gw}{}u\Gd\gr_n dx\right|\leq c\norm u_{L^\infty}\norm{\gr_n}_{W^{2,r'}}\to 0\quad\text{as }n\to\infty,
$$
we get
$$\myint{\Gw}{}u^p dx-\myint{\prt\Gw}{}\myfrac{\prt u}{\prt{\bf n}}dS=M\myint{\Gw}{}|\nabla u|^q dx.
$$
Hence $\nabla u\in L^q(\Gw)$. If $\gz\in C^\infty_0(\Gw)$ and $\gz_n=\gz\eta_n$, there holds
$$-\myint{\Gw}{}\eta_nu\Gd\gz dx+\myint{\Gw}{}\gz u\Gd\gr_n dx+\myint{\Gw}{}u^p\gz_n dx=M\myint{\Gw}{}|\nabla u|^q\gz_n dx.
$$
Since the second term on the left-hand side tends to $0$ and $\gz_n\to\gz$ when $n\to\infty$, we obtain that 
$$-\myint{\Gw}{}u\Gd\gz dx+\myint{\Gw}{}u^p\gz dx=M\myint{\Gw}{}|\nabla u|^q\gz dx.
$$
Hence $u$ is a solution in the sense of distribution in $\Gw$.
\\
Next we show that $\nabla u\in L^2(\Gw)$. Multiplying $(\ref{Z11})$ by $u\eta_n$ and integrating, we obtain
$$\myint{\Gw}{}|\nabla u|^2\eta_n dx-\myint{\Gw}{}u\langle\nabla u,\nabla\gr_n\rangle dx-\myint{\prt\Gw}{}u\myfrac{\prt u}{\prt{\bf n}}dS+\myint{\Gw}{}u^{p+1}\eta_n dx=M\myint{\Gw}{}u|\nabla u|^q\eta_n dx.
$$
As
$$\BA {lll}
\myint{\Gw}{}u\langle\nabla u,\nabla\gr_n\rangle dx=\myfrac{1}{2}\myint{\Gw}{}\langle\nabla u^2,\nabla\gr_n\rangle dx\\[4mm]
\phantom{\myint{\Gw}{}u\langle\nabla u,\nabla\gr_n\rangle dx}
=-\myfrac{1}{2}\myint{\Gw}{} u^2\Gd\gr_n dx,
\EA$$
and
$$\left|\myint{\Gw}{} u^2\Gd\gr_n dx\right|\leq c\norm{u}^2_{L^{\infty}}\norm{\gr_n}_{W^{2,r'}}=o(1)\quad\text{as }\,n\to\infty,
$$
we infer that 
$$\myint{\Gw}{}|\nabla u|^2 dx-\myint{\prt\Gw}{}u\myfrac{\prt u}{\prt{\bf n}}dS+\myint{\Gw}{}u^{p+1} dx=M\myint{\Gw}{}u|\nabla u|^q dx.
$$
Finally if $\gz\in C^\infty_0(\Gw)$ and $\gz_n=\gz\eta_n$, then 
$$
\myint{\Gw}{}\eta_n\langle\nabla u,\nabla\gz\rangle dx-\myint{\Gw}{}\gz\langle\nabla u,\nabla\gr_n\rangle dx
+\myint{\Gw}{}u^p\gz_n dx=M\myint{\Gw}{}|\nabla u|^q\gz_n dx.
$$
Since $r\leq \frac{2N}{N-2}$ there holds
$$\norm{\gr_n}_{W^{1,2}}\leq c\norm{\gr_n}_{W^{2,r'}}\Longrightarrow\norm{\gr_n}_{W^{1,2}}\to 0\quad\text{as }\,n\to\infty.
$$
Using the fact that $\nabla u\in L^2(\Gw)$ and H\"older's inequality, we derive
$$\myint{\Gw}{}\gz\langle\nabla u,\nabla\gr_n\rangle dx\to 0\quad\text{as }\,n\to\infty.
$$
Hence 
$$
\myint{\Gw}{}\langle\nabla u,\nabla\gz\rangle dx
+\myint{\Gw}{}u^p\gz dx=M\myint{\Gw}{}|\nabla u|^q\gz dx.
$$
This implies that $u$ is a weak solution of $(\ref{Z11})$ and it is therefore $C^2$ in $\Gw$. \smallskip

\noindent Next we assume $q=\frac{2p}{p+1}$.    We choose $b=\frac{r-1}{p-1}$ and $(\ref{RW11})$  becomes
  \begin{equation}\label{RW20}
\displaystyle
-\Gd v+\left(1-b-\myfrac{Mp b^{\frac{p-1}{p+1}}\ge^{\frac{p+1}{p}}}{p+1}\right)\myfrac{|\nabla v|^2}{v}+\left(\myfrac{1}{b}-\myfrac{Mb^{\frac{p-1}{p+1}}}{(p+1)\ge^{p+1}}\right)v^{r}\leq 0.
\end{equation}
If $(\ref{RW11*})$ holds with this choice of $b$, $(\ref{RW11**})$ becomes 
  \begin{equation}\label{RW21}\BA {lll}
\displaystyle
1-b-\myfrac{Mp b^{\frac{p-1}{p+1}}\ge^{\frac{p+1}{p}}}{p+1}=1-b-pb\left(\myfrac{M}{p+1}\right)^{\frac{p+1}{p}}\\[0mm]
\phantom{1-b-\myfrac{Mp b^{\frac{p-1}{p+1}}\ge^{\frac{p+1}{p}}}{p+1}}=\myfrac{1}{p-1}\left(p-r-p(r-1)\left(\myfrac{M}{p+1}\right)^{\frac{p+1}{p}}\right).
\EA\end{equation}
If $M<m^{*}_r$ defined by $(\ref{Z12})$, we can choose $\ge$ such that 
$$1-b-\myfrac{Mp b^{\frac{p-1}{p+1}}\ge^{\frac{p+1}{p}}}{p+1}=0,
$$
and
$$\myfrac{1}{b}-\myfrac{Mb^{\frac{p-1}{p+1}}}{(p+1)\ge^{p+1}}=\gt:=\gt(\ge)>0.
$$
Then $v$ satisfies 
$$\BA {lll} -\Gd v+\gt v^r\leq 0\qquad&\text{in }\;\Gw\setminus K.
\EA$$
Since $cap^{\BBR^N}_{ 2,r'}(K)=0$ it follows from \cite{BaPi} that $\displaystyle v\leq \max_{x\in\prt\Gw}v(x)$. Hence $u$ is bounded. The different steps of the proof in the first case applies without any modification: first $\nabla u\in L^q(\Gw)$ and the equation holds in the sense of distributions in $\Gw$, then   $\nabla u\in L^2(\Gw)$ and since $r\leq \frac{2N}{N-2}$ we infer that $u$ is a weak solution and thus a strong one.  $\phantom{-------}$\qeda
%%%%%%%%%%%%%%%%%%%%%%%%%%%%%%%%%%%%%%%%%%%%%%%%%%%%%%%%%%%%%%%%%%%%%%%%%%%%%%%%%%%%%%%%%%%%%%%%%%%%%%%%%%%%%%%%%%%%%%%%%%%%%%%%%%%%%%%%%%%%%%%%%%%%%%%%%%%%%%%%%%%%%%%%%%%%%%%%%%%%%%%%%%%%%%%%%%%%%%%%%%%%%%%%%%%%%%%%%%%%%%%%%%%%%%%%%%%%%%%%%%%%%%%%%%%%%%%%%%%%%%%%%%%%%%%%%%%%

 \mysection{Measure data}
 
 Let $\Gw\subset \BBR^N$ be a bounded smooth domain with diameter smaller than $2R$. 
 Also any Radon measure in $\Gw$ is extended by $0$ in $\Gw^c$ with the same notation. \smallskip

 \subsection{Proof of \rth{meas-1}: the case $1<q<\frac N{N-1}$}
 
 If $1<q<\frac N{N-1}$ assumption $(\ref{Z7})$ with $\gm\geq 0$ reduces to 
 \bel{MD1*}\BA {ll}
\gm(K)\leq Ccap^{{\BBR^N}}_{2,p'}(K)\qquad \text{for all compact set }\,K\subset\Gw.
\EA\ee
 The construction of solutions is based upon the following result due to Boccardo-Murat-Puel \cite{BMP}. 
It is concerned with a general quasilinear equation in a domain $G\subset\BBR^{_N}$
\bel{ChI-3'-3}\BA {ll}
\CQ(u):=-\Gd u+B(.,u,\nabla u)=0\qquad\text{in }\CD'(G),
\EA\ee
where $B\in C(G\ti\BBR\ti\BBR^{_N})$ satisfies
\bel{I-3'-2}\BA {ll}
\abs{B(x,r,\xi)}\leq \Gg(|r|)(1+|\xi|^2)\quad\text{for all }\,(x,r,\xi)\in G\ti\BBR\ti\BBR^{_N},
\EA\ee
for some continuous increasing function $\Gg$ from $\BBR^+$ to $\BBR^+$.

\bth{BMP} Let $G$ be a bounded domain in $\BBR^{_N}$. If there exists a supersolution $\gf$ and a subsolution $\psi$ of the equation $\CQ v=0$
belonging to $W^{1,\infty}(G)$ and such that $\psi\leq\gf$, then for any $\chi\in W^{1,\infty}(G)$ satisfying 
$\psi\leq\chi\leq\gf$ there exists a function $u\in W^{1,2}(G)$ solution of $\CQ u=0$ such that $\psi\leq u\leq\gf$ 
and $u-\chi\in W^{1,2}_0(G)$.\es

The  sub and super solutions are linked to the two problems in which $p$ and $q$ are bigger than $1$, and $\gm$ and $\gw$ are Radon measures
\begin{equation}\label{MD1}\BA{lll}
-\Gd v+|v|^{p-1}v=\gm&\qquad\text{in }\Gw\\
\phantom{-\Gd v+|v|^{p-1}}
v=0&\qquad\text{in }\prt\Gw,
\EA
\end{equation}
and 
\bel{MD2}\BA{lll}
-\Gd w-M|\nabla w|^q=\gw&\qquad\text{in }\Gw\\
\phantom{-\Gd -M|\nabla w|^q}
w=0&\qquad\text{in }\prt\Gw.
\EA
\ee
It is proved in \cite[Theorem 4.1]{ BaPi} that Problem $(\ref{MD1})$ admits a solution, $v\in L^p(\Gw)$, necessarily unique, if and only if $\gm$ is absolutely continuous with respect to the 
Bessel capacity $cap^{\BBR^N}_{2,p'}$, that is 
\bel{MD4}\BA{lll}
\text{\it For any compact set }E\subset\Gw\,, \;cap^{\BBR^N}_{2,p'}(E)=0\Longrightarrow |\gm|(E)=0.
\EA
\ee
Concerning $(\ref{MD2})$, from \cite[Theorem 1.9]{Hu-Ph} a sufficient condition for solvability is the estimate
\bel{MD5}\BA{lll}
\text{\it  For any compact set }E\subset\Gw\,,\; |\gw|(E)\leq Ccap^{\BBR^N}_{1,q'}(E),
\EA
\ee
for some $C>0$. When $\gm$ is nonnegative and has compact support in $\Gw$, this condition turns out to be necessary.
If $(\ref{MD5})$ is satisfied there exists $\ge_0>0$ such that $(\ref{MD2})$ admits a solution with $\gw$ replaced by 
$\ge\gw$ with $0<\ge\leq\ge_0$. Furthermore $\nabla w\in L^q(\Gw)$ and the following estimates hold 
 \cite[Theorem 1.2]{BVHNV},
\bel{MD6}\BA{lll}
|\nabla w(x)|\leq c_{13}\ge {\bf I}^{2R}_1[|\gw|](x),
\EA
\ee
at least if $\gw$ has compact support or is a smooth function, and, with no such conditions on $\gm$,  
\bel{MD6'}\BA{lll}|w(x)|\leq c_{14}\ge {\bf G}^\Gw[|\gw|](x),
\EA
\ee
with $c_{13},c_{14}$ depending on $N$ and $q$, where ${\bf I}^{2R}_1$ is the truncated Riesz potential in $\BBR^N$
 defined for any measure $\gm$ by 
 \bel{riesz1}
{\bf I}^{2R}_1[\gm](x)=\myint{0}{2R}\myfrac{\gm(B_\gr(x))}{\gr^{N-1}}\myfrac{d\gr}{\gr}\quad\text{for all }\;x\in\BBR^N,
 \ee
 and ${\bf G}^\Gw$ the Green potential in $\Gw$. If $R=\infty$ we denote by ${\bf I}_1:={\bf I}^\infty_1$ the classical Riesz potential and if $\Gw=\BBR^N$ the role of ${\bf G}^\Gw$ is played by the Newtownian potential ${\bf I}_2$. We start with the following easy result:
\blemma{cap-eq} Let $r>1$, $k\in \BBN^*$ and $\gm\in\mathfrak M_+(\Gw)$. If $\gm \in W^{-k,r}(\Gw)$ is nonnegative, then there exists
$C>0$ such that
\bel{MD7}\BA{lll}
\gm(E)\leq C\left(cap^{\Gw}_{k,r'}(E)\right)^{\frac {1}{r'}}\quad\text{for any compact set }E\subset\Gw.
\EA
\ee
where $r'=\frac r{r-1}$. Conversely when $k=1,2$ and $\gm$ satisfies 
\bel{MD8}\BA{lll}
\gm(E)\leq Ccap^{\Gw}_{k,r'}(E)\quad\text{for all compact set }E\subset\Gw,
\EA
\ee
for some $C>0$, then, \smallskip

\noindent (i) if $k=2$ then $\gm \in W^{-2,r}(\Gw)$,\smallskip

\noindent (ii) if $k=1$ then $\gm \in W^{-1,r}(\Gw)$.
\es
\Proof Assume first that $\gm\in \mathfrak M_+(\Omega)\cap W^{-k,r}(\Gw)$. Let 
$\gz\in C^\infty_0(\Gw)$ such that $0\leq \gz\leq 1$ and $\gz\geq \chi_{_K}$. Then 
$$\gm(E)\leq \myint{\Gw}{}\gz d\gm=\langle \gm,\gz\rangle\leq \norm\gz_{W^{k,r'}_0}\norm\gm_{W^{-k,r}}.
$$
By the definition of capacity
$$\gm(K)\leq \norm\gm_{W^{-k,r}}\left(cap^{\Gw}_{k,rp'}(K)\right)^{\frac{1}{r'}}.
$$
Conversely if $(\ref{MD8})$ holds with $k=2$ there exists $\ge_0>0$ such that for every $\ge\in (0,\ge_0]$, there exists $z\in L^r(\Gw)$ satisfying
\bel{MD9-}\BA{lll}
-\Gd z=z^r+\ge\gm\qquad&\text{in }\Gw\\
\phantom{-\Gd}
z=0\qquad&\text{on }\prt\Gw,
\EA
\ee
(see \cite[Theorem 2.10, Remark 2.11]{PV}). Since $z\geq {\bf G}^\Gw[\ge\gm]$, it follows that ${\bf G}^\Gw[\gm]\in L^r(\Gw)$ and therefore 
$\gm\in W^{-2,r}(\Gw)$. Since ${\bf G}^\Gw$ is an isomorphism from $L^{r'}(\Gw)$ into $W^{2,r'}(\Gw)\cap W^{1,r'}_0(\Gw)$, we infer by duality that  
${\bf G}^\Gw$ is an isomorphism from $W^{-2,r}(\Gw)$ into $L^{r}(\Gw)$. Hence $\gm \in W^{-2,r}(\Gw)$.\smallskip

Finally, if $(\ref{MD8})$ holds with $k=1$, then there exists $\ge_0>0$ such that for every $\ge\in (0,\ge_0]$ there exists $z\in W^{1,r}(\Gw)$ satisfying
\bel{MD9}\BA{lll}
-\Gd z=|\nabla z|^r+\ge\gm\qquad&\text{in }\Gw\\
\phantom{-\Gd}
z=0\qquad&\text{on }\prt\Gw.
\EA
\ee
Then $z$ satisfies $z\geq {\bf G}^\Gw[\ge\gm]$. Since $z\in L^{r^*}(\Gw)$ by Sobolev imbedding theorem, we have that ${\bf G}^\Gw[\gm]\in L^{r^*}(\Gw)$, which implies the claim. 
\qeda\medskip

\noindent {\it Proof of the theorem}. We put  $\gm_n=\gm\ast\eta_n$ where $\{\eta_n\}\subset C^\infty_0(\BBR^N)$ is a sequence of mollifiers with supp$(\eta_n)\subset B_{\frac 1n}$, and we denote by $v_n$ the solution of 
\bel{MD10}\BA{lll}
-\Gd v+v^p=\ge\gm_n\chi_{_\Gw}\qquad&\text{in }\Gw\\
\phantom{-\Gd +v^p}
v=0\qquad&\text{on }\prt\Gw.
\EA
\ee
Since $\gm$ satisfies $(\ref{MD8})$ so does $\gm_n$ with the same constant $C$. Hence $\gm_n\in W^{-2,p}(\Gw)$ and $\gm_n\to\gm$ in $
W^{-2,p}(\Gw)$ as $n\to \infty$. We also denote by $z_n$ a nonnegative solution of 
\bel{MD11}\BA{lll}
-\Gd z=z^p+\ge\gm_n\chi_{_\Gw}\qquad&\text{in }\Gw\\
\phantom{-\Gd}
z=0\qquad&\text{on }\prt\Gw,
\EA
\ee
and by $w_n$ a nonnegative solution of 
\bel{MD12}\BA{lll}
-\Gd w=M|\nabla w|^q+\ge\gm_n\chi_{_\Gw}\qquad&\text{in }\Gw\\
\phantom{-\Gd}
w=0\qquad&\text{on }\prt\Gw,
\EA
\ee
Since $w_n$ is $C^2$, it is unique by the strong maximum principle. Then there holds by $(\ref{MD6})$, $(\ref{MD6'})$,
\bel{MD13}\BA{lll}
(i)\qquad \qquad&v_n\leq \ge{\bf G}^\Gw[\gm_n]\leq  w_n\leq 
c_{14}\ge{\bf G}^{\Gw}[\gm_n]\leq c_{14}\ge {\bf I}_2[\gm_n]\lfloor_{\Gw},\qquad\qquad\qquad \qquad   \\[2mm]
(ii)\qquad &|\nabla w_n|\leq c_{13}\ge {\bf I}^{2R}_1[\gm_n].
\EA
\ee
Since $v_n$ and $w_n$ are respectively a subsolution and a supersolution of 
\bel{MD14}\BA{lll}
-\Gd u+u^p=M|\nabla u|^q+\ge\gm_n\qquad&\text{in }\Gw\\
\phantom{-\Gd +u^p}
u=0\qquad&\text{on }\prt\Gw,
\EA
\ee
it follows by \rth{BMP} that there exists $u=u_n\in W^{1,\infty}_0(\Gw)$ satisfying $(\ref{MD14})$ in the sense that 
for any $\gz\in C^2_c(\overline\Gw)$ there holds
\bel{MD15}\BA{lll}
-\myint{\Gw}{}u_n\Gd\gz dx+\myint{\Gw}{}\left(u_n^p-M|\nabla u_n|^q\right)\gz dx=\ge\myint{\Gw}{}\gz d\gm_n.
\EA
\ee
It is unique by the strong maximum principle and it satisfies 
\bel{MD16}\BA{lll}
v_n\leq u_n\leq w_n\leq c_{14}\ge {\bf I}_2[\gm_n]\lfloor_{\Gw}.
\EA
\ee
Since ${\bf I}_2[\gm_n]\lfloor_{\Gw}$ is uniformly bounded in $L^p(\Gw)$, the sequence of functions $\{u_n\}$ shares this property.  If $\eta={\bf G}^\Gw[1]$, there holds
\bel{MD16-1}\BA{lll}
\myint{\Gw}{}\left(u_n+\eta u_n^p\right)dx=M\myint{\Gw}{}|\nabla u_n|^q\eta dx+\ge\myint{\Gw}{}\eta d\gm_n.
\EA
\ee
Hence $|\nabla u_n|$ is uniformly bounded in $L^q_{d_{_\Gw}}(\Gw)$ where ${d_{_\Gw}}(x)=\dist (x,\prt\Gw)$. By $(\ref{MD6})$,
\bel{MD16-2}\BA{lll}|\nabla u_n|\leq c_{13}\ge {\bf I}^{2R}_1[|\gm_n-u_n^p|]\leq c_{13}\left(\ge {\bf I}^{2R}_1[\gm_n]+{\bf I}^{2R}_1[u_n^p]\right)\\[2mm]
\phantom{|\nabla u_n|\leq c_{13}\ge {\bf I}^{2R}_1[\gm_n-u_n^p]}
\leq c_{13}\left(\ge {\bf I}^{2R}_1[\gm_n]+c^p_{14}\ge^p{\bf I}^{2R}_1[( {\bf I}_2[\gm_n\lfloor_{\Gw}])^p]\right).
\EA\ee
Using \cite[Lemma 4.2]{HMV} we have equivalence between 
\bel{MD22}
{\bf I}_1[\left({\bf I}_2[\gm_n\lfloor_{\Gw}]\right)^p]\leq c_{15}{\bf I}_1[\gm_n\lfloor_{\Gw}],
\ee
and 
\bel{MD23}
{\bf I}_2[\left({\bf I}_2[\gm_n\lfloor_{\Gw}]\right)^p]\leq c_{17}{\bf I}_2[\gm_n\lfloor_{\Gw}],
\ee
and  $c_{17}\leq c_{15}\leq C(N,p)c_{17}$. Moreover, since diam$(\Gw)<2R$, 
\bel{MD24-1}\left({\bf I}^{2R}_1[\left({\bf I}_2[\gm]\right)^p]\right)^{q}\leq c_{21}^q({\bf I}_1[\gm])^q\leq c^qc_{21}^q({\bf I}^{2R}_1[\gm])^q.
\ee
for some $c=c(N,R)>0$. 
Next, it is quoted in \cite[Theorem 1.1 ]{HMV} that the inequality $(\ref{MD23})$ is equivalent to 
the main assumption of \rth{meas-1},
\bel{MD24}
\gm_n\lfloor_{\Gw}(E)\leq Ccap^{\BBR^N}_{2,p'}(E)\quad\text{for all compact set } E\subset\Gw, 
\ee
for some $C>0$. Actually this equivalence is proved in \cite{MaVer}. By \cite[Theorem 3.14-(a)]{AdHe} there exists $A=A(N)>0$ such that 
\bel{MD16-5}
\left|\{x\in \BBR^N:|{\bf I}_1[\gm_n\lfloor_{\Gw}](x)|>\gl\}\right|\leq A\gl^{-\frac{N}{N-1}}\norm{\gm_n\lfloor_{\Gw}}_{L^1}^{\frac{N}{N-1}}.
\ee
Clearly the above inequality holds if ${\bf I}_1$ is replaced by ${\bf I}_1^{2R}$ and $\BBR^N$ by $\Gw$. This is an estimate of ${\bf I}_1^{2R}[\gm_n\lfloor_{\Gw}]$ in the Lorentz  space $L^{\frac{N}{N-1},\infty}(\Gw)$ (or Marcinkiewicz space). Clearly
$$\norm{\gm_n\lfloor_{\Gw}}_{L^1}\leq \norm{\gm}_{\mathfrak M}.
$$
Therefore $(\ref{MD16-2})$ implies that $|\nabla u_n|$ is bounded in $L^{\frac{N}{N-1},\infty}(\Gw)$, hence equi-integrable in $L^q(\Gw)$ since $q<\frac N{N-1}$.
By \rlemma {cap-eq}-(ii) and classical harmonic analysis results, ${\bf I}_2[\gm_n\lfloor_{\Gw}\to {\bf I}_2[\gm\lfloor_{\Gw}$ in $L^p(\Gw)$ (see e.g. \cite{Ste}). It follows from
 $(\ref{MD16})$ that $u_n$ is equi-integrable in $L^p(\Gw)$. By standard results on elliptic equations and measure theory \cite[Corollary IV]{Br1}, the sequences $\{u_{n}\}$ and $\{|\nabla u_{n}|\}$ are relatively compact in $L^1(\Gw)$. 
Hence there exist a subsequence $\{n_j\}$, converging to $\infty$ and a function $u\in W^{1,q}_0(\Gw)\cap L^p(\Gw)$ such that 
$u_{n_j}\to u$ in $W^{1,1}(\Gw)$ and a.e. in $\Gw$. Since $\{u_{n}\}$ and $\{|\nabla u_{n}|\}$ are also equi-integrable in $L^p(\Gw)$ and $L^q(\Gw)$ respectively, we infer that $u_{n_j}\to u$ in $W^{1,q}_0(\Gw)\cap L^p(\Gw)$. 
 It follows from Vitali's convergence theorem that $u_{n_j}\to u$ in $L^p(\Gw)$ and $\nabla u_{n_j}\to \nabla u$ in $L^q(\Gw)$. Letting 
$n_j\to\infty$ in $(\ref{MD15})$ we conclude that the identity
\bel{MD16-6}\BA{lll}
-\myint{\Gw}{}u\Gd\gz dx+\myint{\Gw}{}\left(u^p-M|\nabla u|^q\right)\gz dx=\ge\myint{\Gw}{}\gz d\gm_n,
\EA
\ee
holds for any $\gz\in C^2_c(\overline\Gw)$. Clearly 
\bel{MD16-7}\BA{lll}
v\leq u\leq w\leq C\ge {\bf I}_2[\gm]\lfloor_{\Gw},
\EA
\ee
where $v$ and $w$ are respectively the solution of $(\ref{MD1})$ and the minimal solution of $(\ref{MD2})$ with $\gm$ replaced by $\ge \gm$.\qeda

 %%%%%%%%%%%%%%%%%%%%%%%%%%%%%%%%%%%%%%%%%%%%%%%%%%%%%%%%%%%%%%%%%%%%%%%%%%%%%%%%%%%%%%%%%%%%%%%%%%%%%%%%%%%%%%%%%%%%%%%%%%%%%%%%%%%%%%%%%%%%%%%%%%%%%%%%%%%%%%%%%%%%%%%%%%%%%%%%%%%%%%%%%%%%%%%%%%%%%%%%%%%%%%
  \subsection{Proof of \rth{meas-1}: the general case}

The approach with super and sub solutions does not work directly and we follow  the method developed for proving \cite[Theorem 1.9]{Hu-Ph} which 
is a delicate extension of the one in the subcritical case. The fact that a sequence of approximation $\{u_n\}$ is bounded in $W^{1,q}_0(\Gw)\cap L^p(\Gw)$ does not imply the uniform integrability of $\{\nabla u_n\}$ in $L^q(\Gw)$ for $q\geq \frac{N}{N-1}$.

\bdef{Mdef} If $r>1$ and $k\in \BBN^*$ we denote by ${\bf M}^{k,r}(\Gw)$ the set of bounded measures $\gm$ in $\Gw$ which satisfy, for some $C>0$,  
\bel{MD9-1}\BA{lll}
|\gm|(E)| \leq Ccap^{\Gw}_{k,r}(E)\quad\text{for all compact sets }E\subset\Gw.
\EA
\ee
If $\Gw$ is replaced by $\BBR^N$, the set is denoted by ${\bf M}^{k,r}$. The smallest constant $C$ such that $(\ref{MD9-1})$ holds is denoted by $[\gm]_{M^{k,r}}$.
\es

For $T>1$ we denote by $E_1(T,\gm)$ the subset of functions $\gz\in W^{1,q}_0(\Gw)$ such that 
\bel{MD17}
\myint{\Gw}{}|\nabla \gz|^qwdx\leq T\ge^q\myint{\Gw}{}({\bf I}^{2R}_1[|\gm|])^q wdx \quad\text{for all }w\in{\bf A}_1\cap L^\infty.
\ee
 
 The following estimate is obtained in \cite[Lemma 5.2]{Hu-Ph} in a more general context.
\blemma{L5.2} Let $q>1$ and $\gm\in {\bf M}^{1,q'}(\Gw)$. Then there exists $c_{17}=c_{17}(N,q)>0$ such that for any $\gz\in E_1(T,\gm)$
\bel{MD18'}
{\bf I}_1^{2R}[|\nabla \gz|^q\chi_{_\Gw}](x)\leq c_{17}T\ge^q[\gm]_{M^{1,q'}}^{q-1}{\bf I}_1^{2R}[|\gm|](x)
\ee
a.e. in $\Gw$.
\es
\blemma{L5.X} Let $q>1$, $\gz\in E_1(T,\gm)$ where $\gm\in {\bf M}^{1,q'}(\Gw)$, and  $S_1(\gz)$ be the solution of
\bel{MD18}
\BA {lll}
-\Gd \phi+\phi^p=M|\nabla \gz|^q+\ge\gm&\qquad\text{in }\Gw\\
\phantom{-\Gd+\phi^p}
\phi=0&\qquad\text{on }\prt\Gw.
\EA
\ee
Then there exists $c_{18}=c_{18}(N,q)>0$ such that 
\bel{MD19}v\leq S_1(\gz)\leq S(\gz)\leq (\ge +c_{18}\ge^qT) I_2[\gm],
\ee
where $S(\gz)$ is the solution of 
\bel{MD20}\BA{lll}
-\Gd \phi=M|\nabla v|^q+\ge\gm\qquad&\text{in }\Gw\\
\phantom{-\Gd}
\phi=0\qquad&\text{on }\prt\Gw,
\EA
\ee
and $v$ is the solution of $(\ref{MD1})$ with $\gm$ replaced by $\ge\gm$. \es
%%%
\Proof There holds
$$S_1(\gz)=\ge {\bf G}^\Gw[\gm]+{\bf G}^\Gw[|\nabla \gz|^q]\leq \ge{\bf I}_2[\gm]+{\bf I}_2[|\nabla \gz|^q]=\ge{\bf I}_2[\gm]+{\bf I}_1[{\bf I}_1[|\nabla \gz|^q]].
$$
By \rlemma{L5.2} with $R=\infty$, ${\bf I}_1[|\nabla \gz|^q\chi_{_\Gw}]\leq c_{17}T\ge^q[\gm]^{q-1}_{M^{1,q'}}{\bf I}_1[\gm]$, hence
\bel{MD21-}\BA{lll}
S_1(\gz)\leq \left(\ge+c_{17}T\ge^q[\gm]^{q-1}_{M^{1,q'}}\right){\bf I}_2[\gm]:=(\ge +c_{18}\ge^qT){\bf I}_2[\gm],
\EA
\ee
since $T\geq 1$ and $\ge\in (0,1]$.\qeda

\blemma{L5.Y}  There exists $T_1>0$ such that for any $T>T_1$, there exists  $\ge_T>0$ such that for all $\ge\in (0,\ge_T]$, $S_1$ maps $E(T,\gm)$ into itself. 
\es
\Proof By \cite[Theorem 1.4]{Hu-Ph}, and using the proof of \cite[Theorem 1.9]{Hu-Ph}, there exist constants $c_{19},c_{20},c_{21}>0$ depending on $N$, $q$, $R$ and $m$ such that for any 
$w\in {\bf A}_1\cap L^\infty$ such that $[w]_{A_1}\leq m$
\bel{MD21}\BA{lll}
\myint{\Gw}{}|\nabla S_1(\gz)|^qwdx\leq c_{19}\myint{\Gw}{}\left({\bf I}^{2R}_1[|\nabla \gz|^q\chi_{_\Gw}]+\ge{\bf I}^{2R}_1[\gm]+{\bf I}^{2R}_1[S^p_1(\gz)]\right)^q w dx\\[4mm]
\phantom{\myint{\Gw}{}|\nabla S_1(\gz)|^qwdx}
\leq c_{20}\myint{\Gw}{}\left(\left({\bf I}^{2R}_1[|\nabla \gz|^q\chi_{_\Gw}]\right)^q+\ge^q({\bf I}^{2R}_1[\gm])^q\right.\\[4mm]
\phantom{---c_{20}\myint{\Gw}{}\left(\left({\bf I}^{2R}_1[|\nabla \gz|^q\chi_{_\Gw}]\right)^q\right)}
\left.+(\ge^{pq}+c_{21}\ge^{pq^2}T^{pq})\left({\bf I}^{2R}_1[\left({\bf I}_2[\gm]\right)^p]\right)^{q}\right) wdx.
\EA
\ee
Using again the equivalence between $(\ref{MD22})$, $(\ref{MD23})$, $(\ref{MD24-1})$ and $(\ref{MD24})$, inequality 
 $(\ref{MD21})$ is transformed in \!
\bel{MD25}\BA{lll}
\myint{\Gw}{}\!|\nabla S_1(\gz)|^qwdx\leq c_{22}\left(\!\ge^q+c_{24}qT^q\ge^{q^2}\!+c_{25}^q(\ge^{pq}\!+c_{20}\ge^{pq^2}T^{pq})\right)\!\myint{\Gw}{} ({\bf I}^{2R}_1[\gm])^qwdx.
\EA
\ee
Finally for any $T>c_{23}$ there exists $\ge_0$ such that for $0<\ge\leq \ge_0$, there exist positive constants $c_{24}$, $c_{25}$ such that 
$$c_{22}\left(\ge^q+c_{24}qT^q\ge^{q^2}+c_{25}^q(\ge^{pq}+c_{20}\ge^{pq^2}T^{pq})\right)\leq T\ge^q,
$$
which implies the claim.\qeda
%%%%
\blemma{L5.Z}  Under the assumptions of \rlemma{L5.Y} with $T>c_{23}$ and $0<\ge\leq\ge_T$, the mapping $S_1$ is compact from $E(T,\gm)$
into itself.
\es
\Proof We prove first the continuity. Let $\{\gz_j\}\subset E(T,\gm)$ such that $\gz_j\to\gz$ in $W^{1,q}_0(\Gw)$.  Using monotonicity as in \cite[Theorem 8]{BrSt} we obtain that for any $r\in [1,\frac N{N-1})$, there exists $\ga=\ga(N,R,r)>0$ such that 
\bel{MD26}\BA{lll}
\ga\norm{S_1(\gz_j)-S_1(\gz_\ell)}_{W^{1,r}_0}+\norm{S^p_1(\gz_j)-S^p_1(\gz_\ell)}_{L^1}\leq M\norm{|\nabla \gz_j|^q-|\nabla \gz_\ell|^q}_{L^1}.
\EA
\ee
Since $\{\nabla \gz_j\}$ is a Cauchy sequence in $L^q(\Gw)$ it follows that $S_1(\gz_j)\to S_1(\gz)$ in $W^{1,r}_0(\Gw)\cap L^p(\Gw)$. Hence 
\bel{MD26^*}\BA{lll}-\Gd S_1(\gz)+S^p_1(\gz)=M|\nabla \gz|^q+\gm&\qquad\text{in }\Gw\\
\phantom{-\Gd +S^p_1(\gz)}S_1(\gz)=0&\qquad\text{in }\prt\Gw.
\EA
\ee
If we take $\gn:=\gn_j=M|\nabla \gz_j|^q+\gm+S_1^p(\gz_j)$. By ($\ref{MD24-1}$) 
$${\bf I}^{2R}_1[S_1^p(\gz_j)]\leq {\bf I}^{2R}_1[\left({\bf I}_2[\gm]\right)^p]\leq C_6{\bf I}^{2R}_1[\gm].
$$
Combined with ($\ref{MD18'}$) we infer that
\bel{MD27'}{\bf I}^{2R}_1[M|\nabla \gz_j|^q+S_1^p(\gz_j)+\gm]\leq c_{28}{\bf I}^{2R}_1[\gm].
\ee  
Let ${\bf M}_1$  be the first order fractional maximal function defined by
\bel{max1}
{\bf M}_1(\gm)(x):=\sup_{\gr>0}\myfrac{\gw(B_\gr(x))}{\gr^{N-1}}\qquad\text{for all }\;x\in\BBR^N.
\ee
 It is classical that for all $\gn\in \mathfrak M_b(\Gw)$, 
\bel{max2}
{\bf M}_1[|\gn|](x)\leq m_{_{^N}}{\bf I}^{2R}_1[|\gn|])(x)\qquad\text{for a.e. }\;x\in\BBR^N.
\ee
Hence, if we set $c_{29}=m_{_N}(c_{28}+1)$, there holds.
\bel{MD27}
{\bf M}_1[M|\nabla \gz_j|^q+\gm+S_1^p(\gz_j)](x)\leq c_{29}{\bf I}^{2R}_1[\gm].
\ee
Since ${\bf I}^{2R}_1[\gm]\in L^q(\Gw)$ we deduce that  the left-hand side of $(\ref{MD27})$ is bounded and equi-integrable in 
$L^q(\Gw)$. Then we apply \cite[Corollary 1.7]{Hu-Ph}, with $w=1$ and deduce that there exists a subsequence $\{\nabla S_1(\gz_{j_n})\}$ which converges in 
$L^q(\Gw)$. By uniqueness the limit is $\{\nabla S_1(\gz)\}$ and the whole sequence $\{\nabla S_1(\gz_{j})\}$ converges. Therefore $S_1$ is continuous. \smallskip

The proof of the compactness follows the same ideas: If $\{\gz_j\}$ is a bounded sequence in $E(T,\gm)$, then $\{\nabla S_1(\gz_j)\}$ is bounded in 
$L^q_w(\Gw)$ by $(\ref{MD21})$, hence $(\ref{MD26})$ holds. By $(\ref{MD23})$, $\{S_1(\gz_j)\}$ is bounded in $L^p(\Gw)$ and equi-integrable since 
${\bf I}_2[\gm]$ belongs to $L^p(\Gw)$ (it is a consequence of $(\ref{MD24})$) and \rlemma{cap-eq}. By \cite{Br1} the sequence 
$\{S_1(\gz_j)\}$ is relatively compact in $W^{1,1}_0(\Gw)$. Hence, there exist $\gz\in E(T,\gm)$ and a subsequence $\{\gz_{j_n}\}$ such that  
$\gz_{j_n}\to \gz$ weakly in $W^{1,q}_0(\Gw)$ and $\{S_1(\gz_{j_n})\}$ converges to some $S\in E(T,\gm)$ in $W^{1,1}_0(\Gw)\cap L^p(\Gw)$, a.e. in $\Gw$ and weakly in $W^{1,q}_0(\Gw)$. As above we derive from \cite[Corollary 1.7]{Hu-Ph} that the sequence $\{{\bf M}_1[M|\nabla \gz_{j_n}|^q+\gm+S_1^p(\gz_j)]\}$ is bounded and equi-integrable in $L^q(\Gw)$, hence $\{\nabla S_1(\gz_{j_n})\}$ converges to $\nabla S$  in $L^q(\Gw)$. Therefore $S$ is a solution of 
\bel{MD28}\BA {lll}
-\Gd S+S^p=M|\nabla \gz|^q+\gm\qquad&\text{in }\Gw\\
\phantom{-\Gd +S^p}
S=0\qquad&\text{on }\prt\Gw.
\EA
\ee
This implies that $S=S_1(\gz)$ and the mapping $S_1$ is compact.\qeda\medskip

 \nind{\it End of the proof of \rth{meas-1}}. It follows from \rlemma{L5.Z} that $S_1$ is a compact mapping from $E(T,\gm)$ into itself. Hence it admits a fixed point $u$
 by Schauder's theorem and $u\in W^{1,q}_0(\Gw)\cap L^p(\Gw)$ is nonnegative and satisfies $(\ref{Z8})$.\qeda
 
 \subsection{Proof of the Corollaries}
 {\it Proof of \rcor{meas-cor1}.} If $q$ satisfies $\frac {Np}{N+p}\leq q<2$, then  $1<p<\frac{2N}{N-2}$. By Sobolev imbedding theorem there holds 
 $$\norm{\gf}_{W^{1,q'}}\leq c_{27}\norm{\gf}_{W^{2,p'}}\quad\text{for all }\gf\in C^2(\overline\Gw),
 $$
where $c_{27}$ depends on $p$, $q$ and $|\Gw|$, provided
$$\myfrac{1}{q'}\geq \myfrac{1}{p'}-\myfrac{1}{N}\Longleftrightarrow q\geq \frac{Np}{N+p}.
$$
The condition $q<2$ necessitates that $\frac{Np}{N+p}<2$, equivalently $p<\frac{2N}{N+2}$. \qeda\medskip

\nind {\it Proof of \rcor{meas-2}.} We recall \cite[Theorem 5.5.1]{AdHe} (a)-(b). If $2p'\leq q'<N$, then
\bel{MD29-}\BA {lll}
\left(cap^{\BBR^N}_{2,p'}(E)\right)^{\frac{1}{N-2p'}}\leq A\left(cap^{\BBR^N}_{1,q'}(E)\right)^{\frac{1}{N-q'}}\qquad\text{for all Borel set }E\subset\Gw.
\EA
\ee
Since $\frac{N-2p'}{N-q'}\geq 1$ we deduce
\bel{MD30}\BA {lll}
cap^{\BBR^N}_{2,p'}(E)\leq A'cap^{\BBR^N}_{1,q'}(E)\qquad\text{for all Borel set }E\subset\Gw.
\EA
\ee
Hence $(\ref{Z9-2})$ implies $(\ref{Z7})$ and existence follows from \rth{meas-1}. The assumption $2p'\leq q'<N$ is equivalent to $\frac N{N-1}<q\leq \frac{2p}{p+1}$. Note that this implies $p>\frac N{N-2}\phantom{}$.$\phantom{------}$\qeda\medskip

\nind {\it Proof of \rcor{meas-3}.} -(i) Since $p$ and $q$ are subcritical, $(\ref{Z7})$ is verified as soon as there holds for every Borel set $E\subset\Gw$, 
$$\gm(E)\leq\myfrac{\gm(\Gw)}{\min\left\{cap_{2,p'}^{\BBR^N}(\{0\}), cap_{1,q'}^{\BBR^N}(\{0\})\right\}}\min\left\{cap_{2,p'}^{\BBR^N}(E), cap_{1,q'}^{\BBR^N}(E)\right\}.
$$
-(ii) If $p$ is subcritical and $q$ supercritical, $(\ref{Z9-1})$ implies $(\ref{Z7})$.\smallskip

\nind -(iii) If $q$ is subcritical and $\gm$ satisfies $(\ref{Z9-2})$ then $(\ref{Z7})$ holds and the existence follows by \rth{meas-1}. However the condition 
that $\gm$ vanishes on Borel sets with $cap_{2,p'}^{\BBR^N}$-capacity is weaker. Let $T>0$ and $B_T$ be the set of $\gz\in W^{1,q}_0(\Gw)$ such that 
$\norm{M|\nabla\gz|^q}_{L^1(\Gw)}\leq T$. If $\gz\in B_T$ there exists a unique solution $S(\gz)\in W^{1,1}_0(\Gw)\cap L^p(\Gw)$ to
\bel{MD30+}\BA {lll}
-\Gd \phi+\phi^p=M|\nabla \gz|^q+\gm\qquad&\text{in }\Gw\\
\phantom{-\Gd +\phi^p}
\phi=0\qquad&\text{on }\prt\Gw.
\EA
\ee
This follows from \cite{BaPi} since the right-hand side of the equation is a measure absolutely continuous with respect to the $cap_{2,p'}^{\BBR^N}$-capacity and the solution is nonnegative. The following estimate holds \cite[Theorem 8]{BrSt}
\bel{MD31}\BA {lll}
(i)\qquad\qquad&\myint{\Gw}{}S^p(\gz)dx\leq \gm(\Gw)+M\myint{\Gw}{}|\nabla \gz|^q dx\qquad\qquad\qquad\qquad\\[4mm]
(ii)\qquad&\ga\norm{S(\gz)}_{W_0^{1,r}}\leq 2\gm(\Gw)+2M\myint{\Gw}{}|\nabla \gz|^q dx,
\EA
\ee
for any $1<r<\frac N{N-1}$. This implies in the case $q=r$, 
$$\left(\myint{\Gw}{}|\nabla S(\gz)|^q dx\right)^{\frac{1}{q}}\leq A\gm(\Gw)+B\myint{\Gw}{}|\nabla \gz|^q dx,
$$
with 
$$A=\frac{2}{\ga}\quad\text{and }\;B=\frac{2M}{\ga}.
$$
For $X>0$ set $F(X)=BX^q-X+A\gm(\Gw)$. Then $F'(X)=0$ iff $X=X_0:=(qB)^{-\frac{1}{q-1}}$ and 
$$F(X_0)=A\gm(\Gw)-\myfrac{q-1}{q}(qB)^{-\frac{1}{q-1}}.
$$
If we assume that 
\bel{MD32}\BA {lll}
F(X_0)< 0\Longleftrightarrow \gm(\Gw)< \myfrac{q-1}{qA}(qB)^{-\frac{1}{q-1}}=(q-1)\left(\myfrac \ga 2\right)^{\frac q{q-1}}M^{-\frac{1}{q-1}},
\EA
\ee
then $\min F(X)<0$, therefore  $F$ admits two positive roots $X_2<X_0<X_1$. Hence the inequality 
\bel{MD33}X_2<\left(\myint{\Gw}{}|\nabla \gz|^q dx\right)^{\frac 1q}\leq X_1,\ee
implies 
\bel{MD34}\left(\myint{\Gw}{}|\nabla S(\gz)|^q dx\right)^{\frac{1}{q}}\leq B\myint{\Gw}{}|\nabla \gz|^q dx+A\gm(\Gw)\leq\left(\myint{\Gw}{}|\nabla \gz|^q dx\right)^{\frac 1q}\leq X_1,
\ee
 and if 
 \bel{MD35}\left(\myint{\Gw}{}|\nabla \gz|^q dx\right)^{\frac 1q}\leq X_2,\ee
then
\bel{MD36}
\left(\myint{\Gw}{}|\nabla S(\gz)|^q dx\right)^{\frac{1}{q}}\leq A\gm(\Gw)+BX^q_2\leq A\gm(\Gw)+BX^q_1=X_1.
\ee
Therefore the mapping $S$ sends $B_{X_1}$ into itself. It is continuous since
\bel{MD37}\BA {lll}
(i)\qquad\quad\quad&\myint{\Gw}{}|S^p(\gz_1)-S^p(\gz_2)|\leq M\myint{\Gw}{}\left||\nabla\gz_1|^q-|\nabla\gz_2|^q\right|dx\\[4mm]
(ii)\qquad\quad\quad&\ga\norm{S(\gz_1)-S(\gz_2)}_{W^{1,r}_0}\leq 2M\myint{\Gw}{}\left||\nabla\gz_1|^q-|\nabla\gz_2|^q\right|dx.
\EA\ee
The fact that this operator is compact follows from the {\it a priori} estimate $(\ref{MD31})$ and \rlemma{L5.Z}. The conclusion is a consequence of Schauder Theorem.
\qeda

\end{document}